\documentclass[final]{siamltex}

\def\qed{\hfill $\diamond$}

\usepackage{amssymb,latexsym,amsmath}
\usepackage{mathrsfs}
\usepackage{enumitem}
\usepackage{caption}
\usepackage{inputenc}
\allowdisplaybreaks[0]

\usepackage{graphics} 
\usepackage{epsfig} 
\usepackage{verbatim}

\newtheorem{assumption}{Assumption}[section]
\newtheorem{remark}{Remark}[section]

\allowdisplaybreaks

\usepackage{color}
\newcommand{\sy}[1]{{\color{black} #1}}
\newcommand{\cm}[1]{{\color{black} #1}}

\usepackage{amsfonts}
\usepackage{float}
\usepackage{multirow}
\usepackage{tikz}
\usetikzlibrary{shapes.geometric, arrows}

\tikzstyle{terminal}=[rectangle, rounded corners, minimum width=1cm, minimum height=1cm,text centered, draw=black]

\tikzstyle{terminal_2}=[rectangle, rounded corners, minimum width=0.5cm, minimum height=0.5cm,text centered, draw=black]

\tikzstyle{terminal_3}=[ellipse, rounded corners, minimum width=0.5cm, minimum height=0.5cm,text centered, draw=black, ultra thick]

\tikzstyle{arrow} = [thick,->,>=stealth]

\tikzstyle{addblock} = [draw,circle]

\tikzstyle{cost_block}=[rectangle, minimum width=1cm, minimum height=1cm,text centered, draw=black, ultra thick]

\begin{document}

\title{Robustness to Incorrect Priors and Controlled Filter Stability in Partially Observed Stochastic Control}

%

\author{Curtis McDonald and Serdar Y\"uksel \thanks{Curtis McDonald is Department of Statistics and Data Science at Yale University, New Haven, CT, USA (email: curtis.mcdonald@yale.edu). S.~Y\"uksel is with the Department of Mathematics and Statistics, Queen's University, Kingston, Ontario, Canada, K7L 3N6 (e-mail: yuksel@mast.queensu.ca). This research was supported in part by the Natural Sciences and Engineering Research Council (NSERC) of Canada. A preliminary version of some of the results was reported at the 2019 IEEE Conference on Decision and Control \cite{mcdonald2019observability}.}}%

\maketitle

\begin{abstract} We study controlled filter stability and its effects on the robustness properties of optimal control policies designed for systems with incorrect priors applied to a true system. Filter stability refers to the correction of an incorrectly initialized filter for a partially observed stochastic dynamical system (controlled or control-free) with increasing measurements. This problem has been studied extensively in the control-free context, and except for the standard machinery for linear Gaussian systems involving the Kalman Filter, few studies exist for the controlled setup. One of the main differences between control-free and controlled partially observed Markov chains is that the filter is always Markovian under the former, whereas under a controlled model the filter process may not be Markovian since the control policy may depend on past measurements in an arbitrary (measurable) fashion. This complicates the dependency structure and therefore results from the control-free literature do not directly apply to the controlled setup. In this paper, we study the filter stability problem for controlled stochastic dynamical systems, and provide sufficient conditions for when a falsely initialized filter merges with the correctly initialized filter over time. These stability results are applied to robust stochastic control problems: under filter stability, we bound the difference in the expected cost incurred for implementing an incorrectly designed control policy compared to an optimal policy. A conclusion is that filter stability leads to much stronger robustness results to incorrect priors (compared with results in \cite{kara2019robustness} without controlled filter stability). \cm{Furthermore, if the optimum cost is that same for each prior, the cost of mismatch between the true \sy{prior} and the assumed \sy{prior} is zero. }
\end{abstract}

\begin{AMS}
93E20, 93E11, 93B07 	
\end{AMS}


\section{Introduction}\label{intro}

In this paper we study partially observed Markov decision problems (POMDP) where a controller/decision maker (DM) has an incorrect prior on the initial state of the system. (i) We first study (controlled) filter stability: the effect this incorrect prior has on the sequence of conditional probabilities on the state variable given the measurements, and when the increasingly available measurement data over time leads the filter to correct itself. (ii) Then we study robustness, when the controller selects a control policy for some cost minimization criterion based on this incorrect prior. We investigate the performance loss due to the mismatch error, that is the difference in expected cost induced by this policy under the true prior compared to that induced by an optimal control policy. 

As we note below, neither of these evidently much related problems has been considered in the literature in the generality considered in this paper. Filter stability has primarily been studied for control-free systems and robustness (to incorrect priors) has been studied without the consideration of filter stability. 

We now describe the probabilistic set up in more detail, see Table \ref{tab_1} for specific notation used to describe random variables, measures, and sigma fields. Let $\mathcal{X}, \mathcal{Y}, \mathcal{U}$ be Polish spaces (that is, complete, separable and metric spaces) equipped with their Borel sigma fields $\mathcal{B}(\mathcal{X}), \mathcal{B}(\mathcal{Y}), \mathcal{B}(\mathcal{U})$; $\mathcal{X}$ will be called the state space, and $\mathcal{Y}$ the measurement space and $\mathcal{U}$ the control action space. We also assume that $\mathcal{X}$ is $\sigma$-compact. 

Define the transition kernel $T$ and observation kernel $Q$ as the mappings 
\begin{table}
\begin{center}
\begin{tabular}{|p{0.4\linewidth}|p{0.6\linewidth}|}
\hline
\textbf{Notation}&\textbf{} \\ \hline
$Y_{[0,n]}$&A finite collection of random variables $\{Y_{0},Y_{1},\cdots,Y_{n}\}$ \\ \hline
$Y_{[0,\infty)}$& An infinite sequence of random variables $Y_{0},Y_{1},\cdots$ \\ \hline
$\mathcal{P}(\mathcal{X})$& The space of probability measures on the space $\mathcal{X}$ with Borel sigma field \sy{$\mathcal{B}(\mathcal{X})$}\\ \hline
$\mu \ll \nu$& The measure $\mu$ is absolutely continuous with respect to $\nu$\\ \hline
$\mathcal{F}^{\mathcal{X}}_{a,b}$& The sigma field generated by $(X_{a},\cdots, X_{b})$\\ \hline
$\mathcal{F}^{\mathcal{X}}_{n}$& The sigma field generated by $X_{n}$\\ \hline
$\mathcal{F}^{\mathcal{X}}_{0,\infty} \vee \mathcal{F}^{\mathcal{Y}}_{0,\infty}$&  The sigma field generated by all state and measurement sequences\\ \hline
$P^{\mu,\gamma}$& The strategic measure on the full sigma field $\mathcal{F}^{\mathcal{X}}_{0,\infty} \vee \mathcal{F}^{\mathcal{Y}}_{0,\infty}$\\ \hline
$P^{\mu,\gamma}(X_{[0,n]}) \equiv P^{\mu,\gamma}|_{\mathcal{F}^{\mathcal{X}}_{0,n}}$& The measure $P^{\mu,\gamma}$ restricted to the sigma field $\mathcal{F}^{\mathcal{X}}_{0,n}$\\ \hline
$P^{\mu,\gamma}((X_{[0,\infty)},Y_{[0,\infty)})|Y_{[0,n]}) \equiv P^{\mu,\gamma}|\mathcal{F}^{\mathcal{Y}}_{0,n}$& The conditional measure of $P^{\mu,\gamma}$ with respect to the sigma field $\mathcal{F}^{\mathcal{Y}}_{0,n}$ \\ \hline
\end{tabular}\end{center}\caption{Some notation for random variables, measures, and sigma fields}\label{tab_1}
\end{table}
\begin{align*}
T:&\mathcal{X}\times \mathcal{U} \to \mathcal{P}(\mathcal{X})
&Q:\mathcal{X} \to \mathcal{P}(\mathcal{Y})\\
&(x,u) \mapsto T(dx'|x,u)& x \mapsto Q(dy|x)
\end{align*}
The POMDP is initialized with a state $X_{0}$ distributed according to a prior $\mu \in \mathcal{P}(\mathcal{X})$. However, the DM does not have access to the state realizations, but instead sees $Y_{0}, Y_1,\cdots, Y_n$ at time $n$. An admissible (deterministic) control policy $\gamma=\{\gamma_{n}\}_{n=0}^{\infty}$ is a sequence of measurable mappings $\gamma_{n}:\mathcal{Y}^{n} \times \mathcal{U}^{n-1} \to \mathcal{U}$ that maps past measurements and control actions to a new control action. At each time stage, the DM uses its control policy to apply a control action $U_{n}=\gamma_{n}(Y_{[0,n]},U_{[0,n-1]})$ which affects the transition kernel to the next state, $X_{n}\sim T(dx|X_{n-1},U_{n-1})$. A new observation is made and the process repeats. Recursively, it follows that $U_{0}$ is a function of $Y_{0}$ and $U_{1}$ is a function of $Y_{0},Y_{1}$, and $U_{0}$. Yet since $U_{0}$ is itself a function of $Y_{0}$, we have that $U_{1}$ is essentially just a function of $Y_{[0,1]}$. In other words, we can restrict ourselves to control policies that are only functions of measurements $Y_{[0,n]}$, though explicit dependence on past actions will be useful for some of our results to follow. \sy{We will denote the collection of {\it admissible} control policies as $\Gamma$}. 

Consider the measurable space $\Omega=\mathcal{X}^{\mathbb{Z}_{+}}\times \mathcal{Y}^{\mathbb{Z}_{+}}$, endowed with the product topology, (that is, $\omega \in \Omega$ is a sequence of states and measurements $\omega=\{(x_{i},y_{i})\}_{i=0}^{\infty}$).
\begin{definition}
For a fixed initial measure $\mu \in \mathcal{P}(\mathcal{X})$ and a policy $\gamma \in \Gamma$, we define the strategic measure $P^{\mu,\gamma}$ as the probability measure on $(\Omega,\mathcal{B}(\Omega))$ such that
\begin{itemize}
\item[i)] For all $A \in \mathcal{B}(\mathcal{X} \times \mathcal{Y})$ we have
\begin{align}
P^{\mu,\gamma}\left((X_{0},Y_{0}) \in A\right)=\int_{A}Q(dy|x)\mu(dx)
\end{align}
\item[ii)] For every $n\geq 1$, for all $A \in \mathcal{B}(\mathcal{X} \times \mathcal{Y})$ let $u_{n-1}=\gamma_{n-1}(y_{[0,n-1]},u_{[0,n-2]})$ then we have
\begin{align}\label{strategic_transition}
P^{\mu,\gamma}((X_{n},Y_{n}) \in A|(X,Y)_{[0,n-1]}=(x,y)_{[0,n-1]}) = \int_{A}Q(dy|x)T(dx|x_{n-1},u_{n-1})
\end{align}
\end{itemize}
\end{definition}
\begin{remark}
Note that $(X,Y)_{[0,\infty)}$ is in general not a Markov chain under \sy{$P^{\mu,\gamma}$} as $u_{n-1}$ depends on the past measurements in equation (\ref{strategic_transition}).
\end{remark}

\subsection{Filter Stability}

Given a prior $\mu\in \mathcal{P}(\mathcal{X})$ and a policy $\gamma \in {\Gamma}$ we can then define the filter and predictor for a POMDP using the strategic measure \sy{$P^{\mu,\gamma}$}.
\begin{definition}
\begin{itemize}
\item[(i)] We define the one step predictor process as the sequence of conditional probability measures
\begin{align*}
\pi_{n-}^{\mu,\gamma}(\cdot)=P^{\mu,\gamma}(X_{n} \in \cdot|Y_{[0,n-1]},U_{[0,n-1]}) = P^{\mu,\gamma}(X_{n} \in \cdot|Y_{[0,n-1]})~~~~n\in \mathbb{N}
\end{align*}
\item[(ii)]We define the filter process as the sequence of conditional probability measures
\begin{align}\label{filterDefn}
\pi_{n}^{\mu,\gamma}(\cdot)=P^{\mu,\gamma}(X_{n} \in \cdot|Y_{[0,n]},U_{[0,n-1]}) =P^{\mu,\gamma}(X_{n} \in \cdot|Y_{[0,n]}), \quad n\in\mathbb{Z}_+
\end{align}
\end{itemize}
\end{definition}
\begin{remark}\label{remark:filter_control}
Recall that the $U_{[0,n-1]}$ are all functions of the $Y_{[0,n-1]}$, so conditioning on the control actions is not necessary in the above definitions. Yet this conditional probability would be {\it policy dependent}; if we condition on the past actions, this conditioning would be {\it policy-independent}. This distinction is important for some of our discussions to follow.
\end{remark} 
\begin{remark}
It will be useful to note that the filter is the strategic measure conditioned on the sigma field $\mathcal{F}^{\mathcal{Y}}_{0,n}$ and restricted to the sigma field $\mathcal{F}^{\mathcal{X}}_{n}$.
\begin{align*}
\pi_{n}^{\mu,\gamma}(\cdot)=P^{\mu,\gamma}(X_{n} \in \cdot|Y_{[0,n]})=P^{\mu,\gamma}|_{\mathcal{F}^{\mathcal{X}}_{n}}|\mathcal{F}^{\mathcal{Y}}_{0,n}
\end{align*}
\end{remark}
Say a prior $\mu \in \mathcal{P}(\mathcal{X})$ and a policy $\gamma \in \Gamma$ are chosen, an observer sees measurements $Y_{[0,\infty)}$ generated via the strategic measure $P^{\mu,\gamma}$. The observer is aware that the policy applied is $\gamma$, but incorrectly thinks the prior is $\nu\neq \mu$. The observer will then compute the incorrectly initialized filter $\pi_{n}^{\nu,\gamma}$ while the true filter is $\pi_{n}^{\mu,\gamma}$. The filter stability problem is concerned with the merging of $\pi_{n}^{\nu,\gamma}$ and $\pi_{n}^{\mu,\gamma}$ as $n$ goes to infinity.

In the literature, there are a number of merging notions when one considers stability which we enumerate here. Let $C_{b}(\mathcal{X})$ represent the set of continuous and bounded functions from $\mathcal{X} \to \mathbb{R}$.
\begin{definition}\label{weak_def}
Two sequences of probability measures $P_{n}$, $Q_{n}$ merge weakly if $\forall~f \in C_{b}(\mathcal{X})$ we have
\begin{small}$\lim_{n \to \infty} \left|\int f dP_{n}-\int f dQ_{n}\right|=0$\end{small}.
\end{definition}
\begin{definition}\label{tv_def}
For two probability measures $P$ and $Q$ we define the total variation norm as 
$\|P-Q\|_{TV}=\sup_{\|f\|_{\infty}\leq 1}\left|\int fdP-\int fdQ \right|$
where $f$ is assumed measurable. We say two sequences of probability measures $P_{n}$, $Q_{n}$ merge in total variation if $\|P_{n}-Q_{n}\|_{TV} \to 0$ as $n \to \infty$.
\end{definition}
\begin{definition}~
\begin{itemize}
\item[(i)] For two probability measures $P$ and $Q$ we define the relative entropy as
$
D(P\|Q)=\int \log \frac{dP}{dQ} dP=\int \frac{dP}{dQ}\log \frac{dP}{dQ}dQ
$
where we assume $P \ll Q$ and $\frac{dP}{dQ}$ denotes the Radon-Nikodym derivative of $P$ with respect to $Q$.
\item[(ii)] Let $X$ and $Y$ be two random variables, let $P$ and $Q$ be two different joint measures for $(X,Y)$ with $P \ll Q$. Then we define the (conditional) relative entropy between $P(X|Y)$ and $Q(X|Y)$ as
\begin{align}
&D(P(X|Y)\|Q(X|Y))=\int  \log \left(\frac{dP_{X|Y}}{dQ_{X|Y}}(x,y)\right)dP(x,y) \nonumber\\
& \qquad \qquad \qquad =\int \left(\int \log \left(\frac{dP_{X|Y}}{dQ_{X|Y}}(x,y)\right)dP(x|Y=y)\right)dP(y)\label{re_cond_def}
\end{align}
\end{itemize}
\end{definition}

Total variation merging implies weak merging, and relative entropy merging \sy{(i.e. $D(P_n\|Q_n) \to 0$)} implies total variation merging via Pinsker's inequality \cite{csiszar1967information}.


\cm{As equation (\ref{re_cond_def}) shows, the relative entropy of two conditional measures is the expectation of the relative entropy of the conditional measures at a fixed conditional value $D(P(X|Y=y)\|Q(X|Y=y))$ where the expectation is taken over the marginal of $P$ on $Y$. In the case of the filters, $\pi_{n}^{\mu, \gamma}$ and $\pi_{n}^{\nu,\gamma}$ play the role of the inner conditional measures, and the strategic measure $P^{\mu,\gamma}$ is the outer marginal measure over $Y_{[0,n]}$. We write this as $E^{\mu,\gamma}[D(\pi_{n}^{\mu,\gamma}\|\pi_{n}^{\nu,\gamma})]$ where $D(\pi_{n}^{\mu,\gamma}\|\pi_{n}^{\nu,\gamma})$ plays the role of the inner integral in (\ref{re_cond_def}).}


\subsection{Control Cost and Robustness}

In addition to filter stability the second, and perhaps more operational, goal of this paper is to utilize stability to study robustness of optimal control to incorrect priors.  We will consider two optimal stochastic control criteria:
\begin{enumerate}
\item The infinite horizon discounted cost
\begin{align*}
J_{\beta}(\mu,\gamma)&=E^{\mu,\gamma}\left[\sum_{i=0}^{\infty}\beta^{i}c(X_{i},U_{i})\right]~~\beta \in [0,1)
\end{align*}
\item The infinite horizon average cost 
\begin{align*}
J_{\infty}(\mu,\gamma)&=\limsup_{N \to \infty}\frac{1}{N}E^{\mu,\gamma}\left[\sum_{i=0}^{N-1}c(X_{i},U_{i}) \right]
\end{align*}
\end{enumerate}

For a given prior $\mu$, we will also consider the optimal costs:
\begin{align*}
J_{\beta}^{*}(\mu)&=\inf_{\gamma \in \Gamma}J_{\beta}(\mu,\gamma), \qquad J_{\infty}^{*}(\mu)=\inf_{\gamma \in \Gamma}J_{\infty}(\mu,\gamma)
\end{align*}

An optimal control policy $\gamma^{\mu}$ for a given prior $\mu$ is the control policy which achieves the infimum of the expected cost over all admissible control policies, 
\begin{align*}
J_{\beta}(\mu,\gamma^{\mu})=\inf_{\gamma \in \Gamma}J_{\beta}(\mu,\gamma) = J^{*}_{\beta}(\mu)
\end{align*}
Consider then if a controller falsely thinks the prior of a system is $\nu$, when in reality the prior is $\mu$. Then the controller will implement the policy $\gamma^{\nu}$ which is optimal with respect to $\nu$ but incurs an expected cost of $J_{\beta}(\mu,\gamma^{\nu})$. If the controller had utilized the correctly designed policy the cost could have been $J_{\beta}^{*}(\mu)$. In studying robustness, we are interested in studying this difference
\begin{align}
J_{\beta}(\mu,\gamma^{\nu})-J_{\beta}^{*}(\mu)\\\label{robust_diff}
J_{\infty}(\mu,\gamma^{\nu})-J_{\infty}^{*}(\mu)
\end{align}
for the discounted cost problem and the average cost problem. We will show that bounds can be derived by utilizing filter stability.

\subsection{Relations Between the Two Problems} 

We now present the connections between filter stability and robustness (in the context of the discounted cost criterion, this also applies to the average cost setup). Consider a prior $\mu$ and a control policy $\gamma^{\nu}$ which is optimal with respect to a different prior $\nu$. For any $n \in \mathbb{N}$, we have
\begin{align}
 &J_{\beta}(\mu,\gamma^{\nu})=E^{\mu,\gamma^{\nu}}[\sum_{i=0}^{\infty}\beta^{i}c(X_{i},U_{i})] \\
 &=E^{\mu,\gamma^{\nu}}[\sum_{i=0}^{n-1}\beta^{i}c(X_{i},U_{i})]+E^{\mu,\gamma^{\nu}}\left[E^{\mu,\gamma^{\nu}}\left[\sum_{i=n}^{\infty}\beta^{i}c(X_{i},U_{i}) |Y_{[0,n-1]}\right] \right] \nonumber\\
&=E^{\mu,\gamma^{\nu}}\left[\sum_{i=0}^{n-1}\beta^{i}c(X_{i},U_{i})\right]  +(\beta^{n})E^{\mu,\gamma^{\nu}}\left[E^{\mu,\gamma^{\nu}}\left[\sum_{i=0}^{\infty}\beta^{i}c(X_{n+i},U_{n+i}) |Y_{[0,n-1]}\right] \right] \nonumber 
\end{align}

\cm{As will be discussed in more detail in the literature review, under mild regularity conditions for a POMDP an optimal control policy $\gamma^{\mu}$ is a stationary function of the filter realization at time $n$, $\pi_{n}^{\mu, \gamma^{\mu}}$. That is, consider some measurable mapping $\Phi:\mathcal{P}(\mathcal{X}) \to \mathcal{U}$ that maps a probability measure to a control action. Recall an optimal control policy $\gamma^{\mu} = \{\gamma^{\mu}\}_{n=0}^{\infty}$ is a sequence of mappings $\gamma^{\mu}:\mathcal{Y}^{N+1} \to \mathcal{U}$. Then there exists some operator $\Phi$ such that each $\gamma^{\mu}_{n}$ mapping is equivalent to this operator acting on the filter realization at time $n$:
\begin{align*}
\gamma^{\mu}(y_{[0,n]})&=\Phi(P^{\mu,\gamma^{\mu}}(X_{n} \in \cdots|Y_{[0,n]} = y_{[0,n]}) = \Phi(\pi_{n}^{\mu,\gamma^{\mu}})
\end{align*}}
\cm{Consider the $0^{th}$ time stage in the problem. $X_{0}$ is distributed according to $\mu$, but the DM thinks the initial prior is $\nu$. An observation $Y_{0}$ is made, and the control action $u_{0}$ is a function of the filter realization believing the prior is $\nu$,
\begin{align*}
u_{0} = \gamma^{\nu}(y_{0}) = \Phi(\pi_{0}^{\nu,\gamma^{\nu}})
\end{align*}
Now consider the $n^{th}$ time stage. Conditioned on $Y_{[0,n-1]}$, $X_{n}$ is distributed according to $\pi_{n-}^{\mu,\gamma^{\nu}}$(the true predictor), but the DM thinks the distribution is the false predictor $\pi^{\nu,\gamma^{\nu}}_{n-}$. An observation $Y_{n}$ is made, and the control action $u_{n}$ is the same stationary function of the filter realization believing the prior is $\nu$
\begin{align*}
u_{n}&=\gamma_{n}^{\nu}(y_{[0,n]})=\Phi(\pi_{n}^{\nu,\gamma^{\nu}})
\end{align*}}
Therefore, the optimal policy under prior $\nu$ at time $n$ is the same as the optimal policy under the prior $\nu'=\pi_{n-}^{\nu,\gamma^{\nu}}$ at time $0$ since it is a stationary function of the filter realization.
\begin{align*}
\gamma^{\nu}_{n}(y_{[0,n]})&=\Phi(\pi_{n}^{\nu,\gamma^{\nu}})=\Phi(\pi_{0}^{v',\gamma^{v'}})
=\gamma^{v'}_{0}(y_{n})
\end{align*}
We see that the true predictor at time $n$, $X_{n}|Y_{[0,n-1]} \sim\pi_{n-}^{\mu,\gamma^{\nu}}$ acts as a new prior for a restarted control problem, and the control policy implemented is optimal with respect to the false predictor $\nu'=\pi_{n-}^{\nu,\gamma^{\nu}}$. We then have
\begin{align*}
E^{\mu,\gamma^{\nu}}\left[\sum_{i=0}^{\infty}\beta^{i}c(X_{i+n},U_{i+n}) |Y_{[0,n-1]}\right]=J_{\beta}(\pi_{n-}^{\mu,\gamma^{\nu}},\gamma^{\nu'})
\end{align*}
and therefore
\begin{align}
J_{\beta}(\mu,\gamma^{\nu})&=E^{\mu,\gamma^{\nu}}\left[\sum_{i=0}^{n-1}\beta^{i}c(X_{i},U_{i})\right]+(\beta^{n})E^{\mu,\gamma^{\nu}}\left[J_{\beta}(\pi_{n-}^{\mu,\gamma^{\nu}},\gamma^{\nu'}) \right] \nonumber
\end{align}
If we instead  apply the correctly designed policy $\gamma^{\mu}$ and let $\mu'=\pi_{n-}^{\mu,\gamma^{\mu}}$ be the correctly initialized predictor we have
\begin{align}
J_{\beta}^{*}(\mu)&=E^{\mu,\gamma^{\mu}}\left[\sum_{i=0}^{n-1}\beta^{i}c(X_{i},U_{i})\right]+(\beta^{n})E^{\mu,\gamma^{\mu}}\left[J_{\beta}(\pi_{n-}^{\mu,\gamma^{\mu}},\gamma^{\mu'}) \right] \nonumber \\
&=E^{\mu,\gamma^{\mu}}\left[\sum_{i=0}^{n-1}\beta^{i}c(X_{i},U_{i})\right]+(\beta^{n})E^{\mu,\gamma^{\mu}}\left[J_{\beta}^{*}(\pi_{n-}^{\mu,\gamma^{\mu}}) \right] \nonumber
\end{align}
The difference satisfies
\begin{align}
 &J_{\beta}(\mu,\gamma^{\nu})-J^{*}_{\beta}(\mu) \nonumber\\
 =& E^{\mu,\gamma^{\nu}}\left[\sum_{i=0}^{n-1}\beta^{i}c(X_{i},U_{i})\right]-E^{\mu,\gamma^{\mu}}\left[\sum_{i=0}^{n-1}\beta^{i}c(X_{i},U_{i})\right] + \beta^{n}\left(E^{\mu,\gamma^{\nu}}\left[J_{\beta}(\pi_{n-}^{\mu,\gamma^{\nu}},\gamma^{\nu'})\right]
 -E^{\mu,\gamma^{\mu}}\left[J^{*}_{\beta}(\pi_{n-}^{\mu,\gamma^{\mu}})\right]\right) \nonumber\\
 =& E^{\mu,\gamma^{\nu}}\left[\sum_{i=0}^{n-1}\beta^{i}c(X_{i},U_{i})\right]-E^{\mu,\gamma^{\mu}}\left[\sum_{i=0}^{n-1}\beta^{i}c(X_{i},U_{i})\right] \nonumber\\
 &\qquad +\beta^{n}\left(E^{\mu,\gamma^{\nu}}\left[J_{\beta}(\pi_{n-}^{\mu,\gamma^{\nu}},\gamma^{\nu'}) + J^{*}_{\beta}(\pi_{n-}^{\mu,\gamma^{\nu}})-J^{*}_{\beta}(\pi_{n-}^{\mu,\gamma^{\nu}})\right] -E^{\mu,\gamma^{\mu}}\left[J^{*}_{\beta}(\pi_{n-}^{\mu,\gamma^{\mu}})\right]\right) \nonumber\\
 =&E^{\mu,\gamma^{\nu}}\left[\sum_{i=0}^{n-1}\beta^{i}c(X_{i},U_{i})\right]-E^{\mu,\gamma^{\mu}}\left[\sum_{i=0}^{n-1}\beta^{i}c(X_{i},U_{i})\right] \label{past_mistakes_cost}\\ 
 &+\beta^{n}\left(E^{\mu,\gamma^{\nu}}\left[J^{*}_{\beta}(\pi_{n-}^{\mu,\gamma^{\nu}})\right]
 -E^{\mu,\gamma^{\mu}}\left[J^{*}_{\beta}(\pi_{n-}^{\mu,\gamma^{\mu}})\right]\right)\label{startegic_measure_cost}\\
 &+\beta^{n}\left(E^{\mu,\gamma^{\nu}}\left[J_{\beta}(\pi_{n-}^{\mu,\gamma^{\nu}},\gamma^{\nu'}) - J^{*}_{\beta}(\pi_{n-}^{\mu,\gamma^{\nu}})\right]\right) \label{approximation_cost}
 \end{align}
 
Therefore, we see that there are in general three costs associated with applying an incorrectly designed control policy to a control system. The first cost (\ref{past_mistakes_cost}) can be thought of as the ``transient'' cost. At time $[0,n-1]$ the control policy $\gamma^{\nu}$ makes different control decisions than the optimal policy $\gamma^{\mu}$. As such, the costs incurred from time $[0,n-1]$ will be different for the optimal control policy and the incorrectly designed policy. The second cost (\ref{startegic_measure_cost}) is the ``strategic measure cost'': $P^{\mu,\gamma^{\nu}}$ and $P^{\mu,\gamma^{\mu}}$ place different expected distributions on $X_{n}$, which acts as the new prior for the system for future stages. Therefore, even if the DM implements optimal control from time $n$ onwards, the ``initial condition'' the DM finds itself in at time $n$ may have a poor optimal cost over the future time stages. The third cost (\ref{approximation_cost}) is the ``filter approximation cost'', which is {\it the term related to filter stability}. Under predictor stability the falsely initialized predictor $\pi_{n-}^{\nu,\gamma^{\nu}}$ and the true predictor $\pi_{n-}^{\mu,\gamma^{\nu}}$ merge as time goes on. These three costs form the fundamentals of studying the losses a DM incurs for using an incorrectly designed policy.

\cm{The rest of the paper is organized as follows. In Section \ref{sec:lit}, we review existing literature related to filter stability, the existence and structure of optimal control policies, and robustness for stochastic control problems. In Section \ref{secMainRes}, we present our main results for filter stability in control systems, as well as our results for robustness derived from filter stability.  Proofs and further discussion of the filter stability results are in Section \ref{ProofFS}, and proofs and discussion of the robustness results are in Section \ref{ProofRobust}. Section \ref{generalizationN} discusses generalizations, some extensions and directions for future research.}

\section{Literature Review}\label{sec:lit}

\subsection{Controlled Filter Stability and Non-Linear Observability} 

Observability is one of the foundational concepts of modern linear systems theory \cite{Caines,goodwin2014adaptive,anderson1979optimal,Kalman,kushner2014partial,KushnerKalmanFilter}. 
For linear systems, exact recovery of any initial condition with measurements available until some finite time is defined as observability and is characterized by an observability rank condition in both continuous and discrete-time \cite{Chen}. For linear systems, such an observability definition is global (as it applies for all initial states) and is universal in the control policies applied, as the control policy does not affect the estimation errors (known as the {\it no-dual effect} \cite{bartse74} property). For non-linear systems, however, due to the challenges in the analysis which prevent globality as well as control-dependence, more modest and localized definitions are to be imposed: For deterministic continuous-time non-linear systems \cite{hermann1977nonlinear} and \cite{sontag1984concept} present local indistinguishability conditions with subtle differences, and establish relations with Lie-theoretic characterizations which generalize observability rank conditions for non-linear systems defined locally. For discrete-time deterministic models, observability has also been defined by invertibility or exact recovery of an initial state, locally, given measurements with finitely many observations. Nijmeier \cite{nijmeijer1982observability} developed discrete-time analogues of the observability notions presented in \cite{hermann1977nonlinear} (see also \cite{sontag1984concept} for sampled continuous-time systems). 

Such definitions, either with local or with exact invertibility, are often too restrictive for stochastic systems driven by noise. Liu and Bitmead \cite{liu2011stochastic,liu2010observability} introduced a non-linear stochastic observability definition through entropy, where the conditional entropy of the hidden state given measurements not being the same as the unconditional entropy implies their observability notion. Ugrinoovski \cite{ugrinovskii2003observability} also presents an information theoretic formulation, and defines observability as an informativeness, but not invertibility, condition, which would lead to a mild notion of observability for non-linear systems.

In the filtering literature for control systems, the classical setup involves the linear Gaussian system. The filter in this case is the celebrated Kalman filter, where the finite-dimensional Kalman filter is computed recursively using the Riccati equation. Under linear observability and controllability conditions, the Riccati equation admits a unique solution \cite{kushner2014partial,Caines}, which is the unique limit of the Riccati recursions regardless of the initialization. Thus, the Kalman Filter is stable with respect to incorrect, though still Gaussian, priors under the aforementioned conditions \sy{(we note that partial convergence and robustness results on the asymptotic equivalence of conditional expectations and linear estimates for non-Gaussian priors for linear systems are reported in \cite{sowers1992discrete})}. 

This concept of being insensitive to incorrect initializations is called {\it filter stability}.



Aside from the above, much of the results on filter stability involves control-free systems. Thus, results have considered partially observed Markov processes (POMP) as opposed to partially observed Markov \textit{decision} processes (POMDP). Since there is no control in such systems, there is no past dependency in the system and the pair $(X_{n},Y_{n})_{n=0}^{\infty}$ is always a Markov chain. For such control-free models, filter stability has been studied extensively and we refer the reader to \cite{chigansky2009intrinsic} 
for a comprehensive review and a collection of different approaches. As discussed in \cite{chigansky2009intrinsic}, filter stability may arise via two separate mechanisms: (i) The transition kernel is in some sense {\it sufficiently} ergodic, forgetting the initial measure and therefore passing this insensitivity (to incorrect initializations) on to the filter process. (ii) The measurement channel is sufficiently informative (or the system is sufficiently observable) about the underlying state. Much of the control-free literature has focused on the first of the two mechanisms noted above. 

For the second type of mechanisms noted above, Chigansky and Liptser \cite{chigansky2006role} and a series of papers by van Handel \cite{van2008discrete,van2009observability,van2009uniform} have presented the first conditions, to our knowledge, where various observability/informativeness conditions lead to filter stability for control-free systems (\sy{see also an earlier study, \cite{stannat2005stability}, for the linear measurements setup}). A related fundamental result which pairs with observability, intrinsically connected with the analysis in the papers noted above, is that of Blackwell and Dubins \cite{blackwell1962merging}: if $P$ and $Q$ are two measures on a fully observed stochastic process $\{X_{n}\}_{n=0}^{\infty}$ with $P \ll Q$, then the conditional distributions on the future based on the past merge in total variation $P$ a.s. , that is $\|P(X_{[n+1,\infty)} \in \cdot|X_{[0,n]})-Q(X_{[n+1,\infty)} \in \cdot|X_{[0,n]})\|_{TV} \to 0 ~~~P~a.s.$, 
this result immediately shows that 
\begin{align*}
\|P^{\mu,\gamma}(Y_{n}|Y_{[0,n-1]})-P^{\nu,\gamma}(Y_{n}|Y_{[0,n-1]})\|_{TV} \to 0 ~~~P^{\mu,\gamma}
\end{align*}
 when $\mu \ll \nu$ and when combined with uniform observability \cite{van2009observability} results in predictor stability. We refer the reader to \cite{chigansky2006role} and \cite{van2009observability} for related results utilizing such a convergence argument though \cite{chigansky2006role} arrive at the convergence result without \cite{blackwell1962merging}. 
 Van Handel also provides a definition of observability for compact state spaces in \cite{van2009observability}: a (control-free) system is observable if every prior results in a unique probability measure on the measurement sequences, $P^{\mu}(Y_{[0,\infty)} \in \cdot)=P^{\nu}(Y_{[0,\infty)} \in \cdot) \implies \mu = \nu$. For non-compact spaces, there are further characterizations \cite{van2009observability}, such as uniform observability.  

Again, for control-free systems, \cite{mcdonald2018stability} introduces an explicit notion of non-linear observability that we continue to work with in this paper. A control-free model is called one step observable if for every $f \in C_{b}(\mathcal{X})$ and every $\epsilon>0$ there exists a measurable and bounded function $g$ such that
\begin{align}
\|f(\cdot)-\int_{\mathcal{Y}}g(y)Q(dy|\cdot)\|_{\infty}<\epsilon
\end{align}
We will show one step observability has an identical application in control systems. \sy{However, for control-free models \cite{mcdonald2018stability} also has a notion of multiple step observability: A POMP is {\it $N$-step observable} if for every 
$f \in C_{b}(\mathcal{X})$ and every $\epsilon>0$ there exists a measurable and bounded function $g$ such that
\begin{align}\label{NstepObs}
\|f(\cdot)-\int_{\mathcal{Y}}g(y_{[1,N]})Q(dy_{[1,N]}|X_{1}=\cdot)\|_{\infty}<\epsilon
\end{align}
A further notion is observability: A POMP is {\it observable} if for every 
$f \in C_{b}(\mathcal{X})$ and every $\epsilon>0$ there exists $N \in \mathbb{N}$ and a measurable and bounded function $g$ such that (\ref{NstepObs}) applies. We will show that such multi-step definitions of observability are more intricate in a controlled setup (POMDP), see Section \ref{generalizedObsContN}.}

Building on a functional analytic duality result between controllability and observability, and along the same spirit as in \cite{van2009observability} with respect to an analysis involving the null space of an appropriate linear map, a recent work on non-linear filter stability is \cite{kim2019lagrangian}.

\cm{We outline the rigorous argument for controlled filter stability in Section \ref{ProofFS}, but here we summarize the main ideas. Blackwell and Dubins \cite{blackwell1962merging} ensures that the conditional measures on the observed chains $\{Y_{n}\}_{n=0}^{\infty}$ merge as time goes to infinity. We would like to trace this back to the unobserved variable, and this is where a notion of observability allows us to conclude the $X_{n}|Y_{[0,n]}$ conditionals merge as well. However, in a control system the past reliance of control policies means $\{X_{n}, Y_{n}\}_{n=0}^{\infty}$ is not a Markov Chain as it is in the uncontrolled case. Nonetheless, the observation channel $Y_{n}|X_{n}\sim Q$ is unaffected by control and is time invariant, thus with our notion of ``one step'' observability we are able to recreate similar merging results in the controlled set up.}


For both control-free and controlled setups, \cite{MYDobrushin2020} studied filter stability through the use of Dobrushin's coefficient involving both of the kernels $T$ and $Q$. Aside from \cite{MYDobrushin2020}, we are not aware of studies which also consider the controlled setup; the following from \cite{MYDobrushin2020} will be utilized in our robustness analysis in the paper (for the discounted cost criterion analysis):

\begin{definition}\cite[Equation 1.16]{dobrushin1956central}
For a kernel operator $K:S_{1} \to \mathcal{P}(S_{2})$ (that is a regular conditional probability from $S_1$ to $S_2$) for standard Borel spaces $S_1, S_2$, we define the Dobrushin coefficient as:
\begin{align}
\delta(K)=\inf\sum_{i=1}^{n}\min(K(x,A_{i}),K(y,A_{i}))\label{Dob_def}
\end{align}
where the infimum is over all $x,y \in S_{1}$ and all partitions $\{A_{i}\}_{i=1}^{n}$ of $S_{2}$.
\end{definition}

\begin{theorem} \cite[Theorem 4.1]{MYDobrushin2020} \label{curtis_result}
Let $\tilde{\delta}(T):=\inf_{u \in \mathcal{U}} \delta(T(\cdot|\cdot,u))$. Assume that for $\mu,\nu \in \mathcal{P}({\cal X})$, we have $\mu\ll\nu$. Then we have
\begin{align*}
E^{\mu,\gamma}\left[\|\pi_{n+1}^{\mu,\gamma}-\pi_{n+1}^{\nu,\gamma}\|_{TV}\right] \leq (1-\tilde{\delta}(T))(2-\delta(Q))E^{\mu,\gamma}\left[\|\pi_{n}^{\mu,\gamma}-\pi_{n}^{\nu,\gamma}\|_{TV}\right].
\end{align*}
In particular, defining $\alpha:=(1-\tilde{\delta}(T))(2-\delta(Q))$, we have
\begin{align*}
E^{\mu,\gamma}\left[\|\pi_{n}^{\mu,\gamma}-\pi_{n}^{\nu,\gamma}\|_{TV}\right]\leq 2\alpha^n.
\end{align*}
\end{theorem}


\subsection{On Regularity of POMDPs and Optimal Policies} It is known that any POMDP can be reduced to a (completely observable) belief-MDP \cite{Yus76}, \cite{Rhe74}, whose states are the posterior state distributions or {\it beliefs} of the observer; that is, the state at time $t$ is the filter realization given in (\ref{filterDefn}).
We call the resulting MDP the belief-MDP\index{Belief-MDP}. If one reduces a POMDP to such a belief-MDP and if one can establish the measurable selection criteria \cite[Chapter 3]{HernandezLermaMCP} (e.g. via showing the weak Feller property), existence of optimal solutions for discounted cost problems would follow. Indeed, recently \cite{FeKaZg14} and \cite{KSYWeakFellerSysCont} have established complementary conditions on when the kernel for the belief-MDP is weakly continuous.

Thus, for the discounted cost criterion, existence of optimal solutions, which are also stationary (in the belief state), follows. Existence also holds for the average cost problem under some further technical conditions facilitating the vanishing discount method \cite{Bor00,Bor07,Bor03,BoBu04}. The convex analytic method of Borkar \cite{borkar2002convex} can also be utilized under the weak continuity conditions presented above (also see \cite{HernandezLermaMCP}); however, for the average cost setup some restrictions may need to be imposed on the initial state \cite{survey}.

Robustness in control refers to the property that a control policy operates well under misspecified system dynamics or models, yet this concept does not have a singular definition in the literature. Robustness problems often involve studies on when a controller has an incorrect system model, or has some uncertainty about the specifics of the actual system model. A common approach in the literature is to design controllers that work with guaranteed performance bounds over a class of uncertain systems under some structured constraints (see \cite{basbern,zhou1996robust}). However, these studies are different from what we consider here in that the controller is {\it aware} that she may have the wrong system specifications and has some limited idea of what the possible models may be. In our problem, we start with two disparate priors $\mu$ and $\nu$ and show that filter stability leads to robustness (and not just continuity as shown in \cite{kara2019robustness}). Thus, our paper provides some unification between controlled filter stability and optimal robust stochastic control. We also note that some applications of filter stability in rigorously approximating optimal control policies by finite memory control policies have been presented in \cite{kara2020near}.

\subsection{Contributions}

\begin{itemize}
\item[(i)] {\bf Controlled filter stability.} We study the filter stability problem for non-linear stochastic systems driven by control. Towards this goal, we present a definition of observability, inspired from a related definition from our prior work for control-free systems. Using this definition, we provide sufficient conditions for when a falsely initialized filter merges with the correctly initialized filter over time. As noted, a primary difference between control-free and  controlled partially observed Markov chains is that, the filter is always Markovian under a control-free model but not so in a controlled setup as the control policy may depend on past measurements in an arbitrary, measurable, manner. We also establish relations between various notions of stability for controlled non-linear filters.

\item[(ii)] {\bf Robustness to incorrect priors under filter stability.} Building on our results on filter stability, we study the robustness problem in terms of the mismatch loss in the expected cost incurred for implementing an incorrectly designed control policy (optimal for an incorrect initialization) compared to an optimal policy (for a correct initialization), and relate filter stability and robustness. In particular, we show that, unlike \cite{kara2019robustness} where only continuity properties of mismatch in the incorrect priors were established, under filter stability even distant incorrect initial priors may lead to correct costs; and under mild conditions, filter stability leads to insensitivity \sy{of mismatch} to incorrect priors, and not just continuity in them.
\end{itemize}


\section{Statements of Main Results}\label{secMainRes}~

We define here the different notions of stability for the filter:

\begin{definition}
\begin{itemize}
\item[(i)] A filter process is said to be stable in the sense of weak merging with respect to a policy $\gamma$ $P^{\mu,\gamma}$ almost surely (a.s.) if there exists a set of measurement sequences $A \subset \mathcal{Y}^{\mathbb{Z}_{+}}$ with $P^{\mu,\gamma}$ probability 1 such that for any sequence in $A$; for any $f \in C_{b}(\mathcal{X})$ and any prior $\nu$ with $\mu \ll \nu$ (i.e., for all Borel $B$ $\nu(B) = 0 \implies \mu(B)=0$) we have
$
\lim_{n \to \infty} \left|\int f d\pi_{n}^{\mu,\gamma}-\int f d\pi_{n}^{\nu,\gamma}\right|=0
$.

\item[(ii)]  A filter process is said to be stable in the sense of total variation in expectation with respect to a policy $\gamma$ if for any measure $\nu$ with $\mu \ll \nu$ we have
$
\lim_{n \to \infty} E^{\mu,\gamma}[\|\pi_{n}^{\mu,\gamma}-\pi_{n}^{\nu,\gamma}\|_{TV}]=0
$.

\item[(iii)]  A filter process is said to be stable in the sense of total variation with respect to a policy $\gamma$ $P^{\mu,\gamma}$ a.s. \cm{if there exists a set of measurement sequences $A \subset \mathcal{Y}^{\mathbb{Z}_{+}}$ with $P^{\mu,\gamma}$ probability 1 such that for any sequence in $A$;} for any measure $\nu$ with $\mu \ll \nu$ we have
$
\lim_{n \to \infty} \|\pi_{n}^{\mu,\gamma}-\pi_{n}^{\nu,\gamma}\|_{TV}=0~~P^{\mu,\gamma}~a.s.
$.

\item[(iv)]  A filter process is said to be stable in the sense of relative entropy with respect to a policy $\gamma$ if for any measure $\nu$ with $\mu \ll \nu$ we have $\lim_{n \to \infty}E^{\mu,\gamma}[D(\pi_{n}^{\mu,\gamma}\|\pi_{n}^{\nu,\gamma})] = 0$.

\item[(v)]  The filter is said to be \textit{universally} stable in one of the above notions if the notion holds with respect to every admissible policy $\gamma \in \Gamma$.

\end{itemize}
\end{definition}

Predictor stability is defined in an analogous fashion for each of the criteria above. These notions of stability are asymptotic notions of stability; they do not imply a rate of convergence. In some of the robustness results to be presented, we need a stronger notion of stability (see Theorem \ref{curtis_result}). 

\begin{definition}\label{exp_stable}
A POMP is said to be universally exponentially stable in total variation if there exists a coefficient $0< \alpha< 1$ such that for any $\mu \ll \nu$ and any policy $\gamma$ we have
\begin{align*}
E^{\mu,\gamma}[\|\pi_{n+1}^{\mu,\gamma}-\pi_{n+1}^{\nu,\gamma}\|_{TV}]\leq \alpha E^{\mu,\gamma}[\|\pi_{n}^{\mu,\gamma}-\pi_{n}^{\nu,\gamma}\|_{TV}]~~n \in \{0,1,\cdots\}
\end{align*}
\end{definition}
\subsection{Controlled Filter Stability Results}\label{CFSSection}
 A key condition that will drive the stability results is the following definition of observability:

\begin{definition}\label{one_step_observability}[A Definition of Observability for Controlled Stochastic Systems]
A POMDP is called one step observable (universal in \sy{admissible} control policies) if for every $f \in C_{b}(\mathcal{X})$ and every $\epsilon>0$ there exists a measurable and bounded function $g$ such that
\begin{align}\label{defnObsOne}
\|f(\cdot)-\int_{\mathcal{Y}}g(y)Q(dy|\cdot)\|_{\infty}<\epsilon
\end{align}
\end{definition}

%
%


\begin{theorem}\label{weak_merging_pred}
Assume that $\mu \ll \nu$ and that the POMDP is one step observable. Then the predictor is universally stable weakly a.s. .
\end{theorem}

In Section \ref{generalizationN}, we present some generalizations but also explain why the proof method used for control-free systems does not apply for controlled systems for more relaxed (in particular, an $N$-step generalization) definitions of observability. A number of examples for measurement channels satisfying Definition \ref{one_step_observability} have been reported in \cite[Section 3]{mcdonald2018stability}.

The observability notion defined above only results in stability of the predictor in the weak sense $P^{\mu,\gamma}$ a.s. .We next extend this stability to total variation $P^{\mu,\gamma}$ a.s. .Let the measurement channel $Q$ be {\it dominated} in the sense that  there exists a reference measure $\lambda$ such that $\forall x \in {\cal X}, Q(Y \in \cdot | x_{n}=x) \ll \lambda(\cdot)$. Then, we define the Radon-Nikodym derivative
\begin{align}\label{dominatingMeasure}
q(x,y):=\frac{dQ(Y_{n} \in \cdot|x_{n}=x)}{d\lambda}(y)
\end{align}
which serves as a likelihood function (a conditional probability density function). We will consider one of the following assumptions.

\begin{assumption}\label{lebesgue_cont_control2}
\begin{itemize}
\item[(i) ]$T(\cdot|x,u)$ is absolutely continuous with respect to a dominating measure $\phi$ for every $x \in \mathcal{X},u \in \mathcal{U}$, so that \cm{$ t(x_1,x,u) =\frac{dT(\cdot|x,u)}{d\phi}(x_{1})$} where $t$ is continuous in $x$ \cm{for every $x_1 \in \mathcal{X}$ and $u \in \mathcal{U}$}.
\item[(ii)] $q(x,y)$ is bounded and continuous in $x$ for every fixed $y$. \sy{Furthermore, $q(x,y) > 0$ for all $x \in \mathcal{X}, y \in \mathcal{Y}$.}
\end{itemize}
\end{assumption}

\begin{assumption}\label{lebesgue_cont_control}
$T(\cdot|x,u)$ is absolutely continuous with respect to a dominating measure $\phi$ for every $x \in \mathcal{X},u \in \mathcal{U}$, so that \cm{$ t(x_1,x,u) =\frac{dT(\cdot|x,u)}{d\phi}(x_{1})$}. The family of (conditional densities) $\{t(\cdot, x,u)\}_{x \in \mathcal{X},u \in \mathcal{U}}$ is uniformly bounded and equicontinuous. 
\end{assumption}

\begin{theorem}\label{pred_tv_thm}
Let $\mu \ll \nu$. Let Assumption \ref{lebesgue_cont_control2} or Assumption \ref{lebesgue_cont_control} hold. If the predictor is universally stable in the weak sense a.s. then it is also universally stable in total variation a.s. .
\end{theorem}

One of the key steps in the proof of Theorem \ref{weak_merging_pred} is that $P^{\mu,\gamma}(Y_{n}\in \cdot|Y_{[0,n-1]})$ and $P^{\nu,\gamma}(Y_{n}\in \cdot|Y_{[0,n-1]})$ merge in total variation $P^{\mu,\gamma}$ a.s. as $n \to \infty$. To achieve this in a POMDP, we apply Blackwell and Dubins \cite{blackwell1962merging} to the measurement process $\{Y_{n}\}_{n=0}^{\infty}$. However, \cite{blackwell1962merging} is fundamentally about predictive measures of the future given the past, and hence only directly implies predictor stability results, not the filter. Filter stability is studied next.


%

\begin{assumption}\label{Ali_assumption}
The measurement channel $Q$ is continuous in total variation. That is, for any sequence $a_{n} \to a \in \mathcal{X}$ we have $\|Q(\cdot|a_{n}) - Q(\cdot|a)\|_{TV} \to 0$ or in other words $ \|P(Y_0 \in \cdot| X_0=a_n)  -  P(Y_0 \in \cdot| X_0=a) \|_{TV} \to 0$.
\end{assumption}

Assumption \ref{lebesgue_cont_control2}(ii), together with the related domination condition (\ref{dominatingMeasure}), implies Assumption \ref{Ali_assumption} (\sy{see \cite[Section 2.3]{KSYWeakFellerSysCont}); see also \cite[Theorem 3]{KSYWeakFellerSysCont} for a partial converse result}.

\begin{theorem}\label{filterstabilityresults}~
\begin{itemize}
\item[(i)] Let Assumption \ref{Ali_assumption} hold. If the predictor is universally stable in weak merging a.s. , then the filter is universally stable in weak merging in expectation. \label{Ali_theorem_control}
\item[(ii)] The filter is universally stable in total variation in expectation if and only if the predictor is universally stable in total variation \sy{in expectation}.
\item[(iii)] The filter is universally stable in total variation in expectation if and only if it is universally stable in total variation  a.s. .
\item[(iv)] Let $\mu \ll \nu$, and assume for any policy $\gamma$ there exists some finite $n$ such that $E^{\mu,\gamma}[D(\pi_{n}^{\mu,\gamma}\|\pi_{n}^{\mu,\gamma})]<\infty$ and some $m$ such that $E^{\mu,\gamma}[D(P^{\mu,\gamma}|_{\mathcal{F}^{\mathcal{Y}}_{0,m}}\|(P^{\nu,\gamma}|_{\mathcal{F}^{\mathcal{Y}}_{0,m}})]<\infty$. Then the filter is universally stable in relative entropy if and only if it is universally stable in total variation in expectation.

\end{itemize}
\end{theorem}

Proofs are given in Section \ref{ProofFS}. 

%
%
%

\subsection{Results on Robustness to Incorrect Priors}\label{robustMainSec}~

Robustness and continuity properties in incorrect priors have been studied in \cite{kara2019robustness} for the single stage cost problem and the infinite horizon discounted cost problem. The paper provides conditions for when robustness difference in equation (\ref{robust_diff}) goes to 0 as $\nu \to \mu$ in either weak convergence or total variation topologies. The paper provides the following useful result:
\begin{theorem}\cite[Theorem 3.2]{kara2019robustness}\label{ali_robust_thm}
Assume the cost function $c$ is bounded, non-negative, and measurable. Let $\gamma^{\nu}$ be the optimal control policy designed with respect to a prior $\nu$. Then we have
\begin{align}\label{aliThm}
\sy{J_{\beta}}(\mu,\gamma^{\nu})-J^{*}_{\beta}(\mu)\leq 2\frac{\|c\|_{\infty}}{1-\beta}\|\mu-\nu\|_{TV}
\end{align}
\end{theorem}
While the average cost case is not studied in \cite{kara2019robustness}, we can adapt the technique to achieve a similar result for the average cost problem (see the appendix for a proof).
\begin{theorem}\label{avg_cost_robust}
Assume the cost function $c$ is bounded, non-negative, and measurable. Let $\gamma^{\nu}$ be the optimal control policy designed with respect to a prior $\nu$. Then we have
\begin{align*}
J_{\infty}(\mu,\gamma^{\nu})-J_{\infty}^{*}(\mu)\leq 2\|c\|_{\infty}\|\mu-\nu\|_{TV}
\end{align*}
\end{theorem}
However, in the following, we will refine these results by considering {\it the effects of filter stability on robustness to incorrect priors}. Unlike \cite{kara2019robustness}, where continuity in priors was studied, here we utilize filter stability to arrive at stronger insensitivity and robustness results.

For the average cost case, asymptotic filter stability (without a rate of convergence) may suffice for robustness, as we establish in the following. 

\begin{theorem}\label{averageRobustness}
Assume the cost function $c$ is bounded, non-negative, and measurable and assume the predictor is universally stable in total variation in expectation. Consider the {\it span} semi-norm $\|\cdot\|_{sp}$:
\begin{align*}
\|J^*_{\infty}\|_{sp}&:=\sup_{\mu_{1} \in \mathcal{P}(\mathcal{X})}J^{*}_{\infty}(\mu_{1})-\inf_{\mu_{2} \in \mathcal{P}(\mathcal{X})}J^{*}_{\infty}(\mu_{2})
\end{align*}
then we have
\begin{align*}
J_{\infty}(\mu,\gamma^{\nu})-J_{\infty}^{*}(\mu)\leq \|J^*_{\infty}\|_{sp}
\end{align*}
\end{theorem}

In particular, if $\|J^*_{\infty}\|_{sp}=0$, then the average cost optimization problem is completely robust to initialization errors (see \cite{Bor00,Bor07,Bor03,BoBu04} for some setups where this holds). 

We next study the discounted cost problem. This result, due to the discounted nature, is not as strong as the one above involving average cost control as it involves exponential filter stability to mitigate the effects of transient costs (which are absent in the average cost criterion setup). 

\begin{theorem}\label{on-line_thm}
Assume the cost function $c$ is bounded, non-negative, and measurable and assume the filter is universally exponentially stable in total variation in expectation with coefficient $\alpha$. Let
\begin{align}
\|J^*_{\beta}\|_{sp}&=\sup_{\mu_{1}\in \mathcal{P}(\mathcal{X})} J_{\beta}^{*}(\mu_{1})-\inf_{\mu_{2} \in \mathcal{P}(\mathcal{X})}J_{\beta}^{*}(\mu_{2}) \nonumber\\
\rho&=\left(\frac{\|c\|_{\infty}}{1-\beta}-\|J^*_{\beta}\|_{sp}\right)\left(\frac{\|c\|_{\infty}}{1-\beta} \right)^{-1}\label{percent_ratio}\\
f(n)&=\beta^{n}(\rho-4\alpha^{n}), \qquad n^{*}=\frac{\ln\left(\left(\frac{\rho}{4}\right)\left(\frac{\ln(\beta)}{\ln(\alpha)+\ln(\beta)}\right) \right)}{\ln \alpha}
\end{align}
For any priors $\mu \ll \nu$, we have
\begin{align}\label{curtisbound}
J_{\beta}(\mu,\gamma^{\nu})-J_{\beta}^{*}(\mu)\leq \frac{\|c\|_{\infty}}{1-\beta}\left(1-\max(f(\lfloor n^{*} \rfloor),f(\lceil n^{*} \rceil) \right)
\end{align}
\end{theorem}

This result may be considered a ``prior independent'' bound since the bound on the robustness distance does not depend on the actual priors $\mu$ and $\nu$ being considered. This is useful when we have no real knowledge of how ``close'' (in total variation distance) $\mu$ and $\nu$ are, thus if they are far apart the bound will not change. 
Comparing (\ref{aliThm}) and (\ref{curtisbound}), we see that when $\|\mu-\nu\|_{TV}$ is small, i.e. our priors start close together to begin with, the bound (\ref{aliThm}) is actually better than our prior independent bound. However, when they are far apart it may be that
\begin{align*}
2\|\mu-\nu\|_{TV}>1-\max(f(\lfloor n^{*} \rfloor),f(\lceil n^{*} \rceil))
\end{align*}
and the prior independent bound is stronger.

\section{Proofs on Predictor and Filter Stability Results}\label{ProofFS}

\subsection{Predictor Stability Results}


Our approach to predictor stability is as follows. By the result of Blackwell and Dubins \cite{blackwell1962merging} we have that the conditional distributions on the measures merge in total variation almost surely for any control policy. Given any measurable and bounded function $g:\mathcal{Y}^{N+1} \to \mathbb{R}$ of multiple measurements we have
\begin{align}
\int g(y_{[n,n+N]})P^{\mu,\gamma}(dy_{[n,n+N]}|Y_{[0,n-1]}) -\int g(y_{[n,n+N]})P^{\nu,\gamma}(dy_{[n,n+N]}|Y_{[0,n-1]}) \label{bd_tv}
\end{align}
goes to zero $P^{\mu,\gamma}$ almost surely as $n \to \infty$ for any control policy $\gamma$. Consider some continuous and bounded  function $f$ on the state space $\mathcal{X}$, in considering the weak merging of the predictor we take the limit of the terms
\begin{align}
|\int f(x_{n})P^{\mu,\gamma}(dx_{n}|y_{[0,n-1]})-\int f(x_{n})P^{\nu,\gamma}(dx_{n}|y_{[0,n-1]})|
\end{align}
Now, if we could replace $f$ with a function $g:\mathcal{Y}^{N+1} \to \mathbb{R}$ integrated over a conditional measure of $P^{\mu,\gamma}(Y_{[0,N]}\in \cdot|X_{0}=x)$, then we would have
\begin{align*}
&\left|\int \left(\int g(\tilde{y}_{[0,N]})P^{\mu,\gamma}(d\tilde{y}_{[0,N]}|X_{0}=\sy{x_{n}})\right)P^{\mu,\gamma}(dx_{n}|y_{[0,n-1]})\right.\\
&\qquad -\left.\int \left(\int g(\tilde{y}_{[0,N]})P^{\mu,\gamma}(d\tilde{y}_{[0,N]}|X_{0}=\sy{x_{n}})\right)P^{\nu,\gamma}(dx_{n}|y_{[0,n-1]}) \right|
\end{align*}
If we can then shown that this is the same as expression (\ref{bd_tv}), then the weak stability of the predictor follows. However, this will require some special properties of the conditional measure of $Y_{[n,n+N]}|X_{n}$. 
Let us consider a level of abstraction higher to better view these properties in more succinct notation. Consider three stochastic processes $A=\{A_{n}\}_{n=0}^{\infty},B=\{B_{n}\}_{n= 0}^{\infty}, C=\{C_{n}\}_{n=0}^{\infty}$ defined on the same measurable space $(\Omega,\mathcal{F})$ mapping to spaces $\mathcal{A},\mathcal{B},\mathcal{C}$ with their respective Borel sigma fields. We can think of $A$ and $C$ as the ``observed'' processes and $B$ as some ``hidden'' process whose realizations are not known by an observer. \cm{For possible measures $P$ and $P^{*}$} on $(\Omega,\mathcal{F})$ we have the following definitions:

\begin{definition}\label{general_conditional_indep}
For a measure $P$ we say the process $C$ only depends on the process $A$ through $B$ if for every $n \in \mathbb{N}$ we have $P$ a.s.
\begin{align*}
&P(C_{n} \in \cdot |A_{n}=a,B_{n}=b)=P(C_{n} \in \cdot |B_{n}=b)
\end{align*}
\end{definition}

\begin{definition}\label{general_time_homogeneous}
For a measure $P$ we say the channel $C|B$ is time homogeneous if for every $n \in \mathbb{N}$, $P$ a.s.
\begin{align*}
P(C_{n} \in \cdot |B_{n}=b)=P(C_{0} \in \cdot |B_{0}=b)
\end{align*}
\end{definition}

\begin{definition}\label{general_measure_equivalent}
For two measures $P$ and $P^{*}$ we say the channel $C|B$ is measure equivalent if for all $n \in \mathbb{N}$ we have, $P,P^{*}$ a.s. :
\begin{align*}
P(C_{n} \in \cdot |B_{n}=b)=P^{*}(C_{n} \in \cdot |B_{n}=b)
\end{align*}
\end{definition}

\begin{definition}\label{general_observable}
For a measure $P$, the channel $C|B$ is observable if for every continuous and bounded function $f:\mathcal{B} \to \mathbb{R}$ \cm{and every  $\epsilon > 0$} we can find a measurable and bounded function $g:\mathcal{C} \to \mathbb{R}$ such that
\begin{align*}
\sup_{b \in \mathcal{B}}\left|f(b)-\int_{\mathcal{C}}g(c_{0})P(dc_{0}|B_{0}=b)\right|<\epsilon
\end{align*}
\end{definition}

\begin{lemma}\label{general_process_lemma}
Let $A$,$B$,$C$ be stochastic processes as above and assume measures $P,P^{*}$ satisfy Definitions \ref{general_conditional_indep}-\ref{general_observable}. Assume that
\begin{align*}
\lim_{n \to \infty}\|P(C_{n}|A_{n})-P^{*}(C_{n}|A_{n})\|_{TV}=0~P~a.s.
\end{align*}
then we have that $P(B_{n}|A_{n})$ and $P^{*}(B_{n}|A_{n})$ merge weakly $P$ a.s. .
\end{lemma}

\textbf{Proof.}
Consider any continuous and bounded function $f:\mathcal{B} \to \mathbb{R}$. Pick any $\epsilon>0$, by observability (Definition \ref{general_observable}) we can find a measurable and bounded function $g:\mathcal{C} \to\mathbb{R}$ such that
\begin{align*}
&\tilde{f}(b)=\int_{\mathcal{C}}g(c_{0})P(dc_{0}|B_{0}=b)&\|f-\tilde{f}\|_{\infty}<\frac{\epsilon}{3}
\end{align*}
now consider
\begin{small}
\begin{align}
&\left|\int f(b_{n})P(db_{n}|a_{n})-\int f(b_{n})P^{*}(db_{n}|a_{n})\right|\nonumber\\
\leq&  \left|\int \tilde{f}(b_{n})P(db_{n}|a_{n})-\int \tilde{f}(b_{n})P^{*}(db_{n}|a_{n}) \right| +\left|\int (f-\tilde{f})(b_{n})P(db_{n}|a_{n}) \right| +\left|\int(f-\tilde{f})(b_{n})P^{*}(db_{n}|a_{n}) \right|\nonumber\\
\leq&   \left|\int \tilde{f}(b_{n})P(db_{n}|a_{n})-\int \tilde{f}(b_{n})P^{*}(db_{n}|a_{n}) \right|+2\|f-\tilde{f}\|_{\infty}\nonumber\\
\leq&  \left|\int \tilde{f}(b_{n})P(db_{n}|a_{n})-\int \tilde{f}(b_{n})P^{*}(db_{n}|a_{n}) \right|+\frac{2}{3}\epsilon \label{general_weak}
\end{align}
\end{small}
we then have
\begin{align*}
&\left|\int_{\mathcal{B}} \tilde{f}(b_{n})P(db_{n}|a_{n})-\int_{\mathcal{B}} \tilde{f}(b_{n})P^{*}(db_{n}|a_{n}) \right|\\
=&\left|\int_{\mathcal{B}} \int_{\mathcal{C}}g(c_{0})P(dc_{0}|B_{0}=b_{n})P(db_{n}|a_{n})  -\int_{\mathcal{B}} \int_{\mathcal{C}}g(c_{0})P(dc_{0}|B_{0}=b_{n})P^{*}(db_{n}|a_{n}) \right|.
\end{align*}
Furthermore, by measure equivalence (Definition \ref{general_measure_equivalent}) we can replace $P(dc_{0}|B_{0}=b_{n})$ with $P^{*}(dc_{0}|B_{0}=b_{n})$ in the second term. By time homogeneity (Definition \ref{general_time_homogeneous}), we can replace $P(dc_{0}|B_{0}=b_{n})$ with $P(dc_{n}|B_{n}=b_{n})$ and the same for $P^{*}$. We then have
\begin{align*}
\bigg|\int_{\mathcal{B}} \int_{\mathcal{C}}g(c_{n})P(dc_{n}|B_{n}=b_{n})P(db_{n}|a_{n}) -\int_{\mathcal{B}} \int_{\mathcal{C}}g(c_{n})P^{*}(dc_{n}|B_{n}=b_{n})P^{*}(db_{n}|a_{n}) \bigg|
\end{align*}
By assumption, $C$ only depends on $A$ through $B$ (Definition \ref{general_conditional_indep}), so we can write $P(dc_{n}|B_{n}=b_{n})=P(dc_{n}|B_{n}=b_{n},A_{n}=a_{n})$ for any $a_{n}$. We finally apply chain rule for conditional probability and we have
\begin{align}
&\left|\int_{\mathcal{C}}g(c_{n})P(dc_{n}|a_{n})-\int_{\mathcal{C}} g(c_{n})P^{*}(dc_{n}|a_{n}) \right|\leq \|P(C_{n}|A_{n})-P^{*}(C_{n}|A_{n})\|_{TV} \label{general_TV}
\end{align}
By assumption, this goes to zero $P$ a.s. as $n\to \infty$ (where the set of convergence of measure $1$ applies for all $g$ by total variation convergence) so we can find an $N$ where $\forall~n>N$ we have (\ref{general_TV}) is less than $\frac{\epsilon}{3}$ and therefore (\ref{general_weak}) is less than $\epsilon$. 
\qed

\textbf{Proof of Theorem \ref{weak_merging_pred}}:
We apply Lemma \ref{general_process_lemma} to our one step observable POMDP. Denote the processes $A_{n}=Y_{[0,n-1]}$, $B_{n}=X_{n}$ and $C_{n}=Y_{n}$. Fix any control policy $\gamma$. We must then check that for the measures $P^{\mu,\gamma}$ and $P^{\nu,\gamma}$ these processes satisfy Definitions \ref{general_conditional_indep} to \ref{general_observable}. The measurement channel $Y_{n}|X_{n}$ is unaffected by the presence of control, and independent of the policy chosen. The distribution of $Y_{n}$ is fully determined by $Q$ and the realization of $X_{n}$.
\begin{align*}
P^{\mu,\gamma}(Y_{n}|Y_{[0,n-1]},X_{n})&=Q(Y_{n}|X_{n})
=P^{\nu,\gamma}(Y_{n}|X_{n},Y_{[0,n-1]})
\end{align*}
therefore the process $Y_{n}$ only depends on $Y_{[0,n-1]}$ through $X_{n}$ (Definition \ref{general_conditional_indep}). Furthermore the channel between $Y_{n}|X_{n}$ is $Q$ regardless of the time index or the initial measure, therefore the channel is time homogeneous (Definition \ref{general_time_homogeneous}) and measure equivalent (Definition \ref{general_measure_equivalent}). The process is observable (Definition \ref{general_observable}) by assumption of being one step observable.
The assumption $\mu \ll \nu$ implies $P^{\mu,\gamma}|_{\mathcal{F}^{\mathcal{Y}}_{0,\infty}} \ll P^{\nu,\gamma}|_{\mathcal{F}^{\mathcal{Y}}_{0,\infty}}$. $\{Y_{n}\}_{n=0}^{\infty}$ is a fully observed stochastic process, therefore by Blackwell and Dubins \cite{blackwell1962merging} we have that
\begin{align}
 \|P^{\mu,\gamma}(Y_{[n,n+N-1]} \in \cdot|Y_{[0,n-1]}) -P^{\nu,\gamma}(Y_{[n,n+N-1]} \in \cdot|Y_{[0,n-1]})\|_{TV} \to 0 \label{BD_zero}
\end{align}
Therefore we satisfy all the conditions of Lemma \ref{general_process_lemma} and the weak merging of the predictor follows. This holds for any $\gamma$, so we have universal stability.
\qed 

We first state the following supporting results.
\begin{lemma}\label{weakPriorweakFilter}
The (measurement-update) map: 
\cm{\begin{align*}
 (\pi^{\mu,\gamma}_{n_-},y_{n}) \mapsto \pi_{n} \quad: \quad \pi^{\mu,\gamma}_{n}(dx_{n}):&=\frac{q(y_{n}|x_{n})\pi^{\mu,\gamma}_{n-}(dx_{n})}{\int_{\mathcal{X}}q(y_{n}|x_{n})\pi_{n-}^{\mu,\gamma}(dx_{n})}
\end{align*}}
which maps from $\mathcal{P}(\mathcal{X}) \times {\cal Y}$ to $\mathcal{P}(\mathcal{X})$ is weakly continuous in $\pi_{n_-}$ for almost every $y$, provided that $q(x,y)$ is \sy{positive}, bounded, and continuous in $x$ for every fixed $y$.
\end{lemma}

\textbf{Proof.} Consider a continuous and bounded $f$ and let $\pi^m_{n_-} \to \pi_{n_-}$ weakly. Then, $P^{\mu,\gamma}$ a.s. ,
\begin{align}
\int f(x_{n}) \frac{q(x_{n},y_{n}) \pi^m_{n_-}(dx_n)} {\int_{{\cal X}} q(x_{n},y_{n}) \pi^m_{n_-}(dx_n)} =  \frac{ \int f(x_{n}) q(x_{n},y_{n}) \pi^m_{n_-}(dx_n)} {\int_{{\cal X}} q(x_{n},y_{n}) \pi^m_{n_-}(dx_n)} \nonumber 
\end{align}
Since $q$ is bounded continuous, both the numerator and the denominator converge. \qed

\begin{lemma}\label{weakPostTVPredictor}
Let $T(dx_1|x,u) = t(x_1,x,u) \phi(dx_1)$ where $t$ is continuous in $x$ for every $x_1$ and $u$. Then, for any policy $\gamma$
the (time-update) map: 
\cm{\[(\pi_{n}^{\mu,\gamma},u_n) \mapsto \pi^{\mu,\gamma}_{n+1-} \quad : \quad  \pi^{\mu,\gamma}_{n+1_-}(dx_{n+1}):=\int_{\mathcal{X}}T(x_{n+1}|x_n,u_n) \pi_{n}^{\mu,\gamma}(dx_n)\] }
which maps from $\mathcal{P}({\cal X}) \times {\cal U}$ to $\mathcal{P}(\mathcal{X})$ is so that
if $\pi^{\nu,\gamma}_n \to \pi^{\mu,\gamma}_n$ weakly then $\pi^{\nu,\gamma}_{n+1-} \to \pi^{\mu,\gamma}_{n+1-}$ in total variation.
\end{lemma}

\textbf{Proof.}
We will build on Scheff\'e's Lemma \cite{Bil86}. For every given history and action, we have
\[\pi^{\nu,\gamma}_{{n+1}_-}(dx_{n+1}) =  \int T(dx_{n+1}|x_n,u_n) \pi^{\nu,\gamma}_{n}(dx_n)\]
Now, $\int T(dx_{n+1}|x_n,u_n)$ is so that for every fixed $u_n$, as \sy{$\pi^{m}_n \to \pi_n$ weakly}
\[ \int_{x_{n+1} \in \cdot} \int_{x_n \in {\cal X}}  t(x_{n+1},x_n,u_n) \phi(dx_{n+1}) \pi^{m}_n(dx_n) \to \int_{x_{n+1} \in \cdot} \int_{x_n \in {\cal X}} t(x_{n+1},x_n,u_n) \phi(dx_{n+1}) \pi_n(dx_n) \]
in total variation since for every fixed $z$, the Radon-Nikodym derivative with respect to $\phi$
\[\frac{d \int  t(x_{n+1},x_n,u_n) \phi(\cdot) \pi^{m}_n(dx_n)} {d\phi}(z) =   \int  t(z,x_n,u_n) \pi^{m}_n(dx_n)\]
satisfies pointwise convergence $\int  t(z,x_n,u_n) \pi^{m}_n(dx_n) \to \int  t(z,x_n,u_n) \pi_n(dx_n)$ (as $\pi^{m}_n \to \pi_n$ weakly), and Scheff\'e's lemma implies that convergence is in total variation. We apply this result to the sequence $\pi^{\nu,\gamma}_n$ converging to $\pi^{\mu,\gamma}_n$.
\qed

\textbf{Proof of Theorem \ref{pred_tv_thm}}

(i) Under Assumption \ref{lebesgue_cont_control2}, the proof follows from Lemma \ref{weakPriorweakFilter} and \ref{weakPostTVPredictor}. \sy{While in Lemma \ref{weakPriorweakFilter} and \ref{weakPostTVPredictor} we consider convergence (and not merging), we note that the proof of Lemma \ref{weakPriorweakFilter} also implies weak merging of the posteriors as the priors weakly merge, and by considering the signed measure $\pi^{\nu,\gamma}_n - \pi^{\mu,\gamma}_n$ in the proof of Lemma \ref{weakPostTVPredictor}, total variation merging is a result of a generalized Scheff\'e's lemma \cite[Theorem 2.8.9]{Bogachev}.}

(ii) Under Assumption \ref{lebesgue_cont_control}, the proof technique follows the approach taken \cite[Theorem 2.10]{mcdonald2018stability} which we will briefly recap and refer the reader for a full consideration. Under Assumption \ref{lebesgue_cont_control}, both predictors $\pi_{n-}^{\mu}$ and $\pi_{n-}^{\nu}$ admit pdf's $f_{n-}^{\mu}$ and $f_{n-}^{\nu}$ with respect to the dominating measure. Define the difference of the pdf's $f_{n-}=f_{n-}^{\mu}-f_{n-}^{\nu}$, which is then $P^{\mu,\gamma}$ a.s. a uniformly bounded and equicontinuous family. We can then apply Arzela-Ascoli theorem \cite{rudin2006real} to show over a compact set, every subsequence of $\{f_{n-}\}_{n=0}^{\infty}$ admits a further subsequence that converges uniformly to 0. Due to the assumed $\sigma$-compactness, the analysis carries over to the entire state space pointwise. Since the measures merge weakly, by \cite[Theorem 8.6.2]{Bog07} total variation merging $P^{\mu,\gamma}$ a.s. follows. \qed

\subsection{Filter Stability Results}

We will present the proof of each of the results in Theorem \ref{filterstabilityresults} separately.

\textbf{Proof of Theorem \ref{filterstabilityresults} i)}:
The proof is the same as \cite[Theorem 2.9]{mcdonald2018stability} with slight notational changes. \qed

%

We now present some supporting results on the Radon Nikodym derivatives. These results are proven in \cite{mcdonald2018stability} for the control-free case (building in part on \cite{Handel}), but carry over with some additional work to the controlled environment. A sketch of the proof of the following lemma is given in the appendix.

\begin{lemma}\label{tv_structure}\cite[Lemma 4.8]{mcdonald2018stability}
Assume $\mu \ll \nu$ and fix any control policy $\gamma$. For any two sigma fields $\mathcal{G}_{1},\mathcal{G}_{2} \subset \mathcal{F}^{\mathcal{X}}_{0,\infty} \vee \mathcal{F}^{\mathcal{Y}}_{0,\infty}$ we have:
\begin{align*}
\|P^{\mu,\gamma}|_{\mathcal{G}_{1}}|\mathcal{G}_{2}-P^{\nu,\gamma}|_{\mathcal{G}_{1}}|\mathcal{G}_{2}\|_{TV}&=\frac{E^{\nu,\gamma}\left[\left.~\left|E^{\nu,\gamma}\left[\frac{d\mu}{d\nu}(X_{0})|\mathcal{G}_{1} \vee \mathcal{G}_{2} \right]-E^{\nu,\gamma}\left[\frac{d\mu}{d\nu}(X_{0})|\mathcal{G}_{2} \right] \right|~\right|\mathcal{G}_{2} \right]}{E^{\nu,\gamma}\left[\frac{d\mu}{d\nu}(X_{0})|\mathcal{G}_{2} \right]}~~~~P^{\mu,\gamma}~a.s.
\end{align*}
\end{lemma}

\begin{theorem}\label{pred_sigma_thm}
Assume $\mu \ll \nu$ and fix any control policy $\gamma$. The predictor is stable in the sense of total variation in expectation with respect to the policy $\gamma$ if and only if
\begin{align}
E^{\nu,\gamma}\left[\left.\frac{d\mu}{d\nu}(X_{0})\right|\bigcap_{n \geq 1}\mathcal{F}_{0,\infty}^{\mathcal{Y}}\vee \mathcal{F}_{n,\infty}^{\mathcal{X}}\right]&=E^{\nu,\gamma}\left[\left.\frac{d\mu}{d\nu}(X_{0})\right|{F}_{0,\infty}^{\mathcal{Y}}\right]~P^{\nu,\gamma}~a.s.
\end{align} 
\end{theorem}
{\bf Proof.}
Begin with the form laid out in Lemma \ref{tv_structure}. We recognize that for the predictor, $\mathcal{G}_{1}=\sigma(X_{n})$ and $\mathcal{G}_{2}=\sigma(Y_{[0,n-1]})$ we have

\begin{align*}
\|\pi_{n-}^{\mu,\gamma}-\pi_{n-}^{\nu,\gamma}\|_{TV}&=\frac{E^{\nu,\gamma}\left[\left.\left|E^{\nu,\gamma}[\frac{d\mu}{d\nu}(X_{0})|Y_{[0,n-1)},X_{n}]-E^{\nu,\gamma}[\frac{d\mu}{d\nu}(X_{0})|Y_{[0,n-1]}]\right|\right|Y_{[0,n-1]}\right]}{E^{\nu,\gamma}\left[\left.\frac{d\mu}{d\nu}(X_{0})\right|Y_{[0,n-1]}\right]}
\end{align*}
Applying Lemma \ref{restrict_structure} leads to
\begin{align*}
&E^{\mu,\gamma}\left[\|\pi_{n-}^{\mu,\gamma}- \pi_{n-}^{\nu,\gamma}\|_{TV} \right]=E^{\nu,\gamma}\left[\frac{dP^{\mu,\gamma}|_{\mathcal{F}^{\mathcal{Y}}_{0,n-1}}}{dP^{\nu,\gamma}|_{\mathcal{F}^{\mathcal{Y}}_{0,n-1}}}\|\pi_{n-}^{\mu,\gamma}- \pi_{n-}^{\nu,\gamma}\|_{TV} \right]\\
&=E^{\nu,\gamma}\left[E^{\nu,\gamma}\left[\left.\frac{d\mu}{d\nu}(X_{0})\right|Y_{[0,n-1]}\right]\|\pi_{n-}^{\mu,\gamma}- \pi_{n-}^{\nu,\gamma}\|_{TV} \right]\\
&=E^{\nu,\gamma}\left[E^{\nu,\gamma}\left[\left.|E^{\nu,\gamma}[\frac{d\mu}{d\nu}(X_{0})|Y_{[0,n-1]},X_{n}]-E^{\nu,\gamma}[\frac{d\mu}{d\nu}(X_{0})|Y_{[0,n-1]}]|\right|Y_{[0,n-1]}\right] \right]\\
&=E^{\nu,\gamma}\left[|E^{\nu,\gamma}[\frac{d\mu}{d\nu}(X_{0})|Y_{[0,n-1]},X_{n}]-E^{\nu,\gamma}[\frac{d\mu}{d\nu}(X_{0})|Y_{[0,n-1]}]|\right]
\end{align*}

We now want to show that
\begin{align}
E^{\nu,\gamma}[\frac{d\mu}{d\nu}(X_{0})|Y_{[0,n-1]},X_{n}]
=E^{\nu,\gamma}[\frac{d\mu}{d\nu}(X_{0})|Y_{[0,\infty]},X_{[n,\infty)}]\label{future_condition}
\end{align}

\begin{figure*}
\begin{center}
\begin{tikzpicture}[scale=0.65, every node/.style={transform shape}]


\node (X0) [terminal_2,  text width = 0.5cm] {$X_{0}$};
\node (X1) [terminal_2, right of = X0, xshift = 2cm, text width = 0.5cm] {$X_{1}$};
\node (X2) [terminal_3, right of = X1, xshift = 2cm, text width = 0.5cm] {$X_{2}$};
\node (X3) [terminal_2, right of = X2, xshift = 2cm, text width = 0.5cm] {$X_{3}$};


\node (Y0) [terminal_3, above of = X0, yshift = 0.5cm, text width = 0.5cm] {$Y_{0}$};
\node (Y1) [terminal_3, above of = X1, yshift = 0.5cm, text width = 0.5cm] {$Y_{1}$};
\node (Y2) [terminal_2, above of = X2, yshift = 0.5cm, text width = 0.5cm] {$Y_{2}$};
\node (Y3) [terminal_2, above of = X3, yshift = 0.5cm, text width = 0.5cm] {$Y_{3}$};

\node (U0) [terminal_3, right of = Y0, xshift = 0.5cm, text width = 0.5cm]{$U_{0}$};
\node (U1) [terminal_3, right of = Y1, xshift = 0.5cm, yshift = 1 cm, text width = 0.5cm]{$U_{1}$};
\node (U2) [terminal_2, right of = Y2, xshift = 0.5cm, yshift = 2 cm, text width = 0.5cm]{$U_{2}$};
\node (U3) [terminal_2, right of = Y3, xshift = 0.5cm, yshift = 3 cm, text width = 0.5cm]{$U_{3}$};

\draw[arrow] (Y0) -- (0,4.5);
\draw[arrow] (Y0) -- (0,3.5);
\draw[arrow] (Y0) -- (0,2.5);
\draw[arrow] (U0) -- (1.5,4.5);
\draw[arrow] (U0) -- (1.5,3.5);
\draw[arrow] (U0) -- (1.5,2.5);
\draw[arrow] (Y1) -- (3,4.5);
\draw[arrow] (Y1) -- (3,3.5);
\draw[arrow] (Y1) -- (3,2.5);
\draw[arrow] (U1) -- (4.5,4.5);
\draw[arrow] (U1) -- (4.5,3.5);
\draw[arrow] (Y2) -- (6,4.5);
\draw[arrow] (Y2) -- (6,3.5);
\draw[arrow] (U2) -- (7.5,4.5);
\draw[arrow] (Y3) -- (9,4.5);

\draw[arrow] (Y0) --  (U0);
\draw[arrow] (0,2.5) --  (U1);
\draw[arrow] (0,3.5) --  (U2);
\draw[arrow] (0,4.5) --  (U3);

\draw[arrow] (X0) -- node[anchor=west] {} (Y0);
\draw[arrow] (X1) -- node[anchor=west] {} (Y1);
\draw[arrow] (X2) -- node[anchor=west] {} (Y2);
\draw[arrow] (X3) -- node[anchor=west] {} (Y3);

\draw[arrow] (X0) -- node[anchor=south, xshift = 1cm] {} (X1);
\draw[arrow] (X1) -- node[anchor=south, xshift = 1cm] {} (X2);
\draw[arrow] (X2) -- node[anchor=south, xshift = 1cm] {} (X3);

\draw[arrow] (U0) -- (1.5,0);
\draw[arrow] (U1) -- (4.5,0);
\draw[arrow] (U2) -- (7.5,0);
\draw[arrow] (U3) -- (10.5,0);

\draw[arrow] (-2 cm, 0 cm) -- node[anchor=south] {$\mu$} (X0);

\draw[arrow] (X3) -- node[anchor=south, xshift = 0.5cm] {$\cdots$} (12,0);

\end{tikzpicture}
\caption{Diagram of Dependence in a POMDP. When $X_{0}$ is conditioned on $X_{2},Y_{[0,1]}$ then $X_{[3,\infty)},Y_{[2,\infty)}$ do not add any new information to the conditioning}
\label{fig:POMDP_dependence}
\end{center}
\end{figure*}
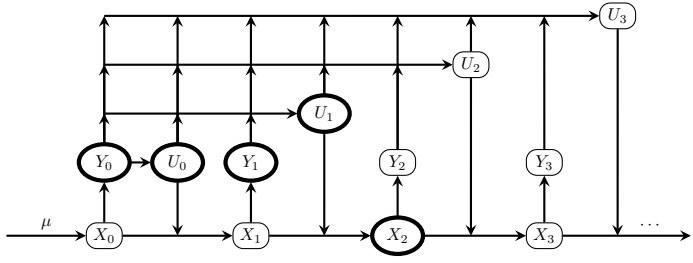

Consider the diagram in Figure \ref{fig:POMDP_dependence} which outlines the dependency structure in a POMDP. Say we are conditioning on $X_{2},Y_{0},Y_{1}$. Since the control policy is known to us, we also know the realizations of $U_{0}$ and $U_{1}$. If we focus on these nodes on the diagram we see that they cordon off $X_{0}$ from $X_{[3,\infty)},Y_{[2,\infty)}$. This is, if we follow back any line from a node $X_{n},n>2$ or $Y_{n},n>1$ to $X_{0}$ we must go through one of these nodes, and hence these future nodes do not add anything useful to the conditioning. Therefore (\ref{future_condition}) holds and we have
\begin{small}
\begin{align}
E^{\mu,\gamma}[\|\pi_{n-}^{\mu,\gamma}-\pi_{n-}^{\nu,\gamma}\|_{TV}]=E^{\nu,\gamma}\left[\left|E^{\nu,\gamma}[\frac{d\mu}{d\nu}(X_{0})|Y_{[0,\infty)},X_{[n,\infty)}]-E^{\nu,\gamma}[\frac{d\mu}{d\nu}(X_{0})|Y_{[0,n-1]}]\right|\right]\label{martingale_expression}
\end{align}
\end{small}

We then see that $A_{n}=E^{\nu}[\frac{d\mu}{d\nu}(X_{0})|Y_{[0,n-1]}]$ is a non-negative uniformly integrable martingale (with respect to the measure $P^{\nu}$) adapted to the increasing filtration $\mathcal{F}^{\mathcal{Y}}_{0,n-1}$. Hence the limit as $n \to \infty$ in $L^{1}(P^{\nu})$ is $E^{\nu}[\frac{d\mu}{d\nu}(X_{0})|\mathcal{F}^{\mathcal{Y}}_{0,\infty}]$. Similarly, we can view $B_{n}=E^{\nu}[\frac{d\mu}{d\nu}(X_{0})|Y_{[0,\infty)},X_{[n,\infty)}]$ as a backwards non-negative uniformly integrable martingale (with respect to the measure $P^{\nu}$) adapted to the decreasing sequence of filtrations $\mathcal{F}_{0,\infty}^{\mathcal{Y}}\vee \mathcal{F}_{n,\infty}^{\mathcal{X}}$. Then by the backwards martingale convergence theorem, the limit as $n \to \infty$ in $L^{1}(P^{\nu})$ is $E^{\nu}[\frac{d\mu}{d\nu}(X_{0})|\bigcap_{n=1}^{\infty}\mathcal{F}^{\mathcal{Y}}_{\cm{0},\infty}\vee \mathcal{F}^{\mathcal{X}}_{n,\infty}]$. \qed

\begin{theorem}\label{filt_tv_thm}
Assume $\mu \ll \nu$ and fix any control policy $\gamma$. The filter is stable in the sense of total variation in expectation with respect to the policy $\gamma$ if and only if
\begin{align}
E^{\nu,\gamma}\left[\left.\frac{d\mu}{d\nu}(X_{0})\right|\bigcap_{n \geq 0}\mathcal{F}_{0,\infty}^{\mathcal{Y}}\vee \mathcal{F}_{n,\infty}^{\mathcal{X}}\right]&=E^{\nu,\gamma}\left[\left.\frac{d\mu}{d\nu}(X_{0})\right|{F}_{0,\infty}^{\mathcal{Y}}\right]~P^{\nu,\gamma}~a.s.
\end{align} 
\end{theorem}

{\bf Proof of Theorem \ref{filterstabilityresults} ii)}
The sigma fields $\mathcal{F}^{\mathcal{X}}_{n,\infty} \vee \mathcal{F}^{\mathcal{Y}}_{0,\infty}$ are a decreasing sequence, that is $\mathcal{F}^{\mathcal{X}}_{n+1,\infty} \vee \mathcal{F}^{\mathcal{Y}}_{0,\infty}~\subset~\mathcal{F}^{\mathcal{X}}_{n,\infty} \vee \mathcal{F}^{\mathcal{Y}}_{0,\infty}$. Therefore, when we take their intersection, removing the first or largest sigma field $\mathcal{F}^{\mathcal{X}}_{0,\infty} \vee \mathcal{F}^{\mathcal{Y}}_{0,\infty}$ from the intersection of a decreasing set of sigma fields does not change the overall intersection. Fix any control policy $\gamma$. The conditions in Theorem \ref{pred_sigma_thm} and \ref{filt_tv_thm} are the same for each policy $\gamma$, therefore if either the filter or the predictor is universally stable in total variation in expectation then the other is as well. \qed 
\qed

{\bf Proof of Theorem \ref{filterstabilityresults} iii)}
Via Fatou's Lemma for the limit supremum, it is clear that almost sure stability implies in expectation. We focus on the other direction.

Assume the filter is universally stable in total variation in expectation. Consider the measure $\rho=\frac{\mu+\nu}{2}$. We have that $\mu \ll \rho$ and $\nu \ll \rho$ and $\|\frac{d\mu}{d\rho}\|_{\infty}\leq 2,\|\frac{d\nu}{d\rho}\|_{\infty}\leq 2$. By Lemma \ref{tv_structure} and the analysis in the proof of Theorem \ref{pred_sigma_thm} \cm{but applied to the filter},
\begin{align*}
\|\pi_{n}^{\mu,\gamma}-\pi_{n}^{\rho,\gamma}\|_{TV}&=\frac{E^{\rho,\gamma}\left[\left.\left|E^{\rho,\gamma}[\frac{d\mu}{d\rho}(X_{0})|Y_{[0,\infty)},X_{[n,\infty)}]-E^{\rho,\gamma}[\frac{d\mu}{d\rho}(X_{0})|Y_{[0,n]}]\right|\right|Y_{[0,n]}\right]}{E^{\rho,\gamma}\left[\left.\frac{d\mu}{d\rho}(X_{0})\right|Y_{[0,n]}\right]}\\
\|\pi_{n}^{\nu,\gamma}-\pi_{n}^{\rho,\gamma}\|_{TV}&=\frac{E^{\rho,\gamma}\left[\left.\left|E^{\rho,\gamma}[\frac{d\nu}{d\rho}(X_{0})|Y_{[0,\infty)},X_{[n,\infty)}]-E^{\rho,\gamma}[\frac{d\nu}{d\rho}(X_{0})|Y_{[0,n]}]\right|\right|Y_{[0,n]}\right]}{E^{\rho,\gamma}\left[\left.\frac{d\nu}{d\rho}(X_{0})\right|Y_{[0,n]}\right]}
\end{align*}

Since the Radon Nikodym derivatives are finite, the expressions in the denominators and numerators above are uniformly integrable martingales with respect to the measure $P^{\rho,\gamma}$ and hence admit limits $P^{\rho,\gamma}$ a.s. . This implies that they admit limits $P^{\mu,\gamma}$ a.s. . Therefore $\|\pi_{n}^{\mu,\gamma}-\pi_{n}^{\rho,\gamma}\|_{TV}$ and $\|\pi_{n}^{\cm{\nu},\gamma}-\pi_{n}^{\cm{\rho},\gamma}\|_{TV}$ admit limits $P^{\mu,\gamma}$ a.s. . By dominated convergence theorem,
\begin{align*}
E^{\mu}[\lim_{n \to \infty} \|\pi_{n}^{\mu,\gamma}-\pi_{n}^{\rho,\gamma}\|_{TV}]&= \lim_{n \to \infty}E^{\mu}[ \|\pi_{n}^{\mu.\gamma}-\pi_{n}^{\rho,\gamma}\|_{TV}]=0
\end{align*}
since we assume the filter is universally stable in expectation by assumption. Total variation is always non-negative, therefore $\lim_{n \to \infty} \|\pi_{n}^{\mu,\gamma}-\pi_{n}^{\rho,\gamma}\|_{TV}$ must be zero a.s. . Similarly,
\begin{align*}
& E^{\mu}[\lim_{n \to \infty} \|\pi_{n}^{\nu,\gamma}-\pi_{n}^{\rho,\gamma}\|_{TV}]=\lim_{n \to \infty}E^{\mu}[ \|\pi_{n}^{\nu,\gamma}-\pi_{n}^{\rho,\gamma}\|_{TV}]\\
&\qquad \qquad \leq \lim_{n \to \infty}E^{\mu}[ \|\pi_{n}^{\nu,\gamma}-\pi_{n}^{\mu,\gamma}\|_{TV}]+\lim_{n \to \infty}E^{\mu}[ \|\pi_{n}^{\mu,\gamma}-\pi_{n}^{\rho,\gamma}\|_{TV}] =0 
\end{align*}
so $\lim_{n \to \infty} \|\pi_{n}^{\nu,\gamma}-\pi_{n}^{\rho,\gamma}\|_{TV}$ must be 0 a.s. . Then we have
\begin{align*}
& \limsup_{n \to \infty} \|\pi_{n}^{\mu,\gamma}-\pi_{n}^{\nu,\gamma}\|_{TV}\leq \limsup_{n \to \infty}  \|\pi_{n}^{\mu,\gamma}-\pi_{n}^{\rho,\gamma}\|_{TV} +\|\pi_{n}^{\mu,\gamma}-\pi_{n}^{\nu,\gamma}\|_{TV}\\
&\qquad \qquad \leq \limsup_{n \to \infty}  \|\pi_{n}^{\mu,\gamma}-\pi_{n}^{\rho,\gamma}\|_{TV} +\limsup_{n \to \infty}\|\pi_{n}^{\mu,\gamma}-\pi_{n}^{\nu,\gamma}\|_{TV} =0
\end{align*}
since the limit supremum of the sum  of two non-negative sequences is less than the sum of the limit supremum of the sequences. Therefore, the limit supremum of $\|\pi_{n}^{\mu}-\pi_{n}^{\nu}\|_{TV}=0$ $P^{\mu,\gamma}$ a.s. . Thus the limit exists and is 0. \qed 


{\bf Proof of Theorem \ref{filterstabilityresults} iv)} First assume the filter is universally stable in relative entropy and fix some policy $\gamma$. Since the square root function is continuous and convex, we have
\begin{align*}
0=\lim_{n \to \infty} \sqrt{\frac
{2}{\log(e)}E^{\mu,\gamma}[D(\pi^{\mu,\gamma}_{n}\|D(\pi_{n}^{\nu,\gamma})]}\geq
\lim_{n \to \infty} E^{\mu,\gamma}\left[\sqrt{\frac
{2}{\log(e)}D(\pi^{\mu,\gamma}_{n}\|D(\pi_{n}^{\nu,\gamma})} \right]
\end{align*}
where we have applied Jensen's inequality. We then apply Pinsker's inequality and we have $\lim_{n \to \infty}E^{\mu,\gamma}[\|\pi_{n}^{\mu,\gamma}-\pi_{n}^{\nu,\gamma}\|_{TV}]=0$.

For the converse direction; assume the filter is universally stable in total variation in expectation and fix a policy $\gamma$. If we apply the chain rule for relative entropy, we have

\begin{align*}
E^{\mu,\gamma}[D(\pi_{n}^{\mu,\gamma}\|\pi_{n}^{\nu,\gamma})]&=D(P^{\mu,\gamma}(X_{n}|Y_{[0,n]})\|P^{\nu,\gamma}(X_{n}|Y_{[0,n]})))\\
&=D(P^{\mu,\gamma}(X_{n},Y_{[0,n]})\|P^{\nu,\gamma}(X_{n},Y_{[0,n]}))-D(P^{\mu,\gamma}(Y_{[0,n]})\|P^{\nu,\gamma}(Y_{[0,n]}))
\end{align*}
We also have by chain rule
\begin{align*}
&D(P^{\mu,\gamma}(X_{n},Y_{[0,n]})\|P^{\nu,\gamma}(X_{n},Y_{[0,n]}))=D(P^{\mu,\gamma}(X_{[n,\infty)},Y_{[0,\infty]})\|P^{\nu,\gamma}(X_{[n,\infty)},Y_{[0,\infty)}))\\
& \qquad \qquad \qquad \qquad -D(P^{\mu,\gamma}((X,Y)_{[n+1,\infty)} | X_{n},Y_{[0,n]}) \sy{\|} P^{\sy{\nu},\gamma}((X,Y)_{[n+1,\infty)}|X_{n},Y_{[0,n]}))
\end{align*}
Consider now Figure \ref{fig:POMDP_dependence}. If we fix $X_{2}, Y_{[0,2]}$ then we also know the control actions $U_{[0,2]}$ then the distribution of $X_{3}$ is fully determined by $X_{2},U_{2}$ and the transition kernel, it is independent of the prior measure. Then the distribution of $Y_{3}$ is determined by $X_{3}$ and $Q$, and $U_{3}$ and $X_{4}$ and so on. Therefore, the conditional measures $P^{\mu,\gamma}((X,Y)_{[n+1,\infty)}|X_{n},Y_{[0,n]})$ and \sy{$P^{\nu,\gamma}((X,Y)_{[n+1,\infty)}|X_{n},Y_{[0,n]})$} are the same even though the priors are different, and therefore the relative entropy is zero.

We then have
\begin{align*}
E^{\mu,\gamma}[D(\pi_{n}^{\mu,\gamma} \| \pi_{n}^{\nu,\gamma})]
&=D(P^{\mu,\gamma}(X_{[n,\infty)},Y_{[0,\infty)})\|P^{\nu,\gamma}(X_{[n,\infty)},Y_{[0,\infty)}))-D(P^{\mu,\gamma}(Y_{[0,n]})\|P^{\nu,\gamma}(Y_{[0,n]}))\\
&=D(P^{\mu,\gamma}|_{\mathcal{F}^{\mathcal{X}}_{n,\infty} \vee \mathcal{F}^{\mathcal{Y}}_{0,\infty}}\|P^{\nu,\gamma}|_{\mathcal{F}^{\mathcal{X}}_{n,\infty} \vee \mathcal{F}^{\mathcal{Y}}_{0,\infty}})-D(P^{\mu,\gamma}|_{\mathcal{F}^{\mathcal{Y}}_{[0,n]}}\|P^{\nu,\gamma}|_{\mathcal{F}^{\mathcal{Y}}_{[0,n]}})
\end{align*}
For each of these, building on \cite[Theorem 3]{BarronSorrento}, we have that the convergence to zero follows the arguments presented in \cite[Theorem 2.12]{mcdonald2018stability}.  \qed
\section{Proofs on Robustness Results}\label{ProofRobust}~

\textbf{Proof of Theorem \ref{averageRobustness}}
Fix some finite $n$, we have
\begin{align*}
& J_{\infty}(\mu,\gamma^{\nu})=\limsup_{T \to \infty}\frac{1}{T}\left(\sum_{i=0}^{n-1}E^{\mu,\gamma^{\nu}}[c(X_{i},U_{i})]+\sum_{i=n}^{T-1}E^{\mu,\gamma^{\nu}}[c(X_{i},U_{i})]\right)\\
&\leq \limsup_{T \to \infty}\frac{1}{T}\sum_{i=0}^{n-1}E^{\mu,\gamma^{\nu}}[c(X_{i},U_{i})]+\limsup_{T \to \infty}\frac{1}{T}E^{\mu,\gamma^{\nu}}\left[\sum_{i=n}^{T-1}c(X_{i},U_{i})\right]\\
&\leq \limsup_{T \to \infty} \frac{n\|c\|_{\infty}}{T}+\limsup_{T \to \infty}\frac{1}{T}E^{\mu,\gamma^{\nu}}\left[\sum_{i=0}^{T-n-1}c(X_{n+i},U_{n+i}) \right]\\
&=\limsup_{T \to \infty}\frac{1}{T}E^{\mu,\gamma^{\nu}}\left[\sum_{i=0}^{T-n-1}c(X_{n+i},U_{n+i}) \right]
\end{align*}
therefore, we see that no matter what decision the DM makes in the first $n$ time stages, since $n$ is finite and $c$ bounded this cost will eventually be dominated by the denominator as $T \to \infty$ and there will be no transient cost associated with this robustness problem. We then claim that
\begin{align}
\limsup_{T \to \infty}\frac{1}{T}E^{\mu,\gamma^{\nu}}\left[\sum_{i=0}^{T-n-1}c(X_{n+i},U_{n+i}) \right]&=\limsup_{T \to \infty}\frac{1}{T-n}E^{\mu,\gamma^{\nu}}\left[\sum_{i=0}^{T-n-1}c(X_{n+i},U_{n+i}) \right] \nonumber 
\end{align}
All terms in the two limsup expressions are positive and bounded since $c$ is a non-negative bounded function, therefore we have in the following that (i) the difference of the limsups is less than or equal to the limsup of the difference, and (ii) the limsup of a product is less than or equal to the product of the limsups. Using these results, 
\begin{align*}
&\limsup_{T \to \infty}\frac{1}{T-n}E^{\mu,\gamma^{\nu}}\left[\sum_{i=0}^{T-n-1}c(X_{n+i},U_{n+i}) \right]-\limsup_{T \to \infty}\frac{1}{T}E^{\mu,\gamma^{\nu}}\left[\sum_{i=0}^{T-n-1}c(X_{n+i},U_{n+i}) \right]\\
&\leq \limsup_{T \to \infty}\left(\frac{1}{T-n}-\frac{1}{T}\right)E^{\mu,\gamma^{\nu}}\left[\sum_{i=0}^{T-n-1}c(X_{n+i},U_{n+i}) \right]\\
& \leq \left(\limsup_{T \to \infty}\frac{n}{T} \right)\left(\limsup_{T \to \infty}\frac{1}{T-n}E^{\mu,\gamma^{\nu}}\left[\sum_{i=0}^{T-n-1}c(X_{n+i},U_{n+i}) \right] \right) =0
\end{align*}
We then apply iterated expectations and Fatou's lemma and we have
\begin{align*}
&\limsup_{T \to \infty}\frac{1}{T-n}E^{\mu,\gamma^{\nu}}\left[\sum_{i=0}^{T-n-1}c(X_{n+i},U_{n+i}) \right] = \limsup_{T \to \infty}\frac{1}{T-n}E^{\mu,\gamma^{\nu}}\left[E^{\mu,\gamma^{\nu}}\left[\sum_{i=0}^{T-n-1}c(X_{n+i},U_{n+i})|Y_{[0,n-1]}\right] \right]\\
\leq &E^{\mu,\gamma^{\nu}}\left[\limsup_{T \to \infty}\frac{1}{T-n}E^{\mu,\gamma^{\nu}}\left[\sum_{i=0}^{T-n-1}c(X_{n+i},U_{n+i})|Y_{[0,n-1]}\right] \right]\
\end{align*}
The optimal control policy is a time invariant function of the filter realization so the predictor at time $n$ acts as a new prior for the control problem. We then have
\begin{align*}
E^{\mu,\gamma^{\nu}}\left[\limsup_{T \to \infty}\frac{1}{T-n}E^{\mu,\gamma^{\nu}}\left[\sum_{i=0}^{T-n-1}c(X_{n+i},U_{n+i})|Y_{[0,n-1]}\right] \right]&=E^{\mu,\gamma^{\nu}}\left[J_{\infty}(\pi_{n-}^{\mu,\gamma^{\nu}},\gamma^{\pi_{n-}^{\nu,\gamma^{\nu}}}) \right]
\end{align*}
With this established, we can now move on to our robustness problem. 
\begin{align}
&J_{\infty}(\mu,\gamma^{\nu})-J_{\infty}^{*}(\mu)\leq E^{\mu,\gamma^{\nu}}[J_{\infty}(\pi_{n-}^{\mu,\gamma^{\nu}},\gamma^{\pi_{n-}^{\nu,\gamma^{\nu}}})]-\inf_{\tilde{\mu} \in \mathcal{P}(\mathcal{X})} J_{\infty}^{*}(\tilde{\mu}) \nonumber \\
&=E^{\mu,\gamma^{\nu}}[J_{\infty}(\pi_{n-}^{\mu,\gamma^{\nu}},\gamma^{\pi_{n-}^{\nu,\gamma^{\nu}}})+J_{\infty}^{*}(\pi_{n-}^{\mu,\gamma^{\nu}})-J_{\infty}^{*}(\pi_{n-}^{\mu,\gamma^{\nu}})]-\inf_{\tilde{\mu} \in \mathcal{P}(\mathcal{X})} J_{\infty}^{*}(\tilde{\mu}) \nonumber \\
&=E^{\mu,\gamma^{\nu}}[J_{\infty}(\pi_{n-}^{\mu,\gamma^{\mu}},\gamma^{\pi_{n-}^{\nu,\gamma^{\nu}}})-J_{\infty}^{*}(\pi_{n-}^{\mu,\gamma^{\nu}})]+E^{\mu,\gamma^{\nu}}[J_{\infty}^{*}(\pi_{n-}^{\mu,\gamma^{\nu}})]-\inf_{\tilde{\mu} \in \mathcal{P}(\mathcal{X})} J_{\infty}^{*}(\tilde{\mu})\\
&\leq E^{\mu,\gamma^{\nu}}[J_{\infty}(\pi_{n-}^{\mu,\gamma^{\nu}},\gamma^{\pi_{n-}^{\nu,\gamma^{\nu}}})-J_{\infty}^{*}(\pi_{n-}^{\mu,\gamma^{\nu}})]+\sup_{\tilde{\mu} \in \mathcal{P}(\mathcal{X})}J_{\infty}^{*}(\tilde{\mu})-\inf_{\tilde{\mu} \in \mathcal{P}(\mathcal{X})} J_{\infty}^{*}(\tilde{\mu}) \nonumber \\
&=E^{\mu,\gamma^{\nu}}[J_{\infty}(\pi_{n-}^{\mu,\gamma^{\nu}},\gamma^{\pi_{n-}^{\nu,\gamma^{\nu}}})-J_{\infty}^{*}(\pi_{n-}^{\mu,\gamma^{\nu}})]+\|J^*_{\infty}\|_{sp} \label{keyArgumentRobust}
\end{align}
by Theorem \ref{avg_cost_robust} we have
\begin{align}
J_{\infty}(\mu,\gamma^{\nu})-J_{\infty}^{*}(\mu)\leq 2\|c\|_{\infty}E^{\mu,\gamma^{\nu}}[\|\pi_{n-}^{\mu,\gamma^{\nu}}-\pi_{n-}^{\nu,\gamma^{\nu}}\|_{TV}]+\|J^*_{\infty}\|_{sp}\label{n_bound}
\end{align}
and this result holds for any $n$ since our choice of $n$ was arbitrary. By assumption, the predictor is universally stable in total variation in expectation. Therefore, for any $\epsilon>0$ there exists an $N$ such that for all $n>N$ we have $E^{\mu,\gamma^{\nu}}[\|\pi_{n-}^{\mu,\gamma^{\nu}}-\pi_{n-}^{\nu,\gamma^{\nu}}\|_{TV}]\leq \frac{\epsilon}{2\|c\|_{\infty}}$
and therefore since (\ref{n_bound}) holds for every $n$, $J_{\infty}(\mu,\gamma^{\nu})-J_{\infty}^{*}(\mu)\leq \|J^*_{\infty}\|_{sp}+\epsilon$
for any $\epsilon>0$, yielding our result.
\qed

\textbf{Proof of Theorem \ref{on-line_thm}}
Pick any $n \in \mathbb{N}$. Starting from expressions (\ref{past_mistakes_cost}), (\ref{startegic_measure_cost}), and (\ref{approximation_cost}) we will consider the three costs. The transient cost is upper bound by
\begin{align*}
E^{\mu,\gamma^{\nu}}\left[\sum_{i=0}^{n-1}\beta^{i}c(x_{i},u_{i})\right]-E^{\mu,\gamma^{\mu}}\left[\sum_{i=0}^{n-1}\beta^{i}c(x_{i},u_{i})\right] \leq \|c\|_{\infty}\sum_{i=0}^{n-1}\beta^{i}=\|c\|_{\infty}\left(\frac{1-\beta^{n}}{1-\beta} \right)
\end{align*}
the strategic measure cost is upper bound by
\begin{align*}
\beta^{n}\left(E^{\mu,\gamma^{\nu}}\left[J^{*}_{\beta}(\pi_{n-}^{\mu,\gamma^{\nu}})\right]
 -E^{\mu,\gamma^{\mu}}\left[J^{*}_{\beta}(\pi_{n-}^{\mu,\gamma^{\mu}})\right]\right)\leq \beta^{n}\|J^*_{\beta}\|_{sp}
\end{align*}
and the approximation cost satisfies 
\begin{align*}
\beta^{n}\left(E^{\mu,\gamma^{\nu}}\left[J_{\beta}(\pi_{n-}^{\mu,\gamma^{\nu}},\gamma^{\pi_{n-}^{\nu,\gamma^{\nu}}}) - J^{*}_{\beta}(\pi_{n-}^{\mu,\gamma^{\nu}})\right]\right) \leq 4 \frac{\|c\|_{\infty}}{1-\beta}(\alpha \beta)^{n}
\end{align*}
Putting these together
\begin{align*}
J_{\beta}(\mu,\gamma^{\nu})-J_{\beta}^{*}(\mu)&\leq \|c\|_{\infty}\left(\frac{1-\beta^{n}}{1-\beta} \right)+\beta^{n}\|J^*_{\beta}\|_{sp}+4 \frac{\|c\|_{\infty}}{1-\beta}(\alpha \beta)^{n}\\
&=\frac{\|c\|_{\infty}}{1-\beta}+\beta^{n} \left(\|J^*_{\beta}\|_{sp}+4\frac{\|c\|_{\infty}}{1-\beta}\alpha^{n}-\frac{\|c\|_{\infty}}{1-\beta} \right)\\
&=\frac{\|c\|_{\infty}}{1-\beta}\left(1+\beta^{n}(4\alpha^{n}-\rho) )\right) =\frac{\|c\|_{\infty}}{1-\beta}(1-f(n))
\end{align*}
This holds for any $n$. Taking the derivative in $n$, at $n^*$, it follows that this will be maximized at some $n^* \in \mathbb{R}_+$, 
\begin{align*}
\alpha^{n^{*}}=\left(\frac{\rho}{4}\right)\left(\frac{\ln(\beta)}{\ln(\alpha)+\ln(\beta)}\right), \qquad n^{*}=\frac{\ln\left(\left(\frac{\rho}{4}\right)\left(\frac{\ln(\beta)}{\ln(\alpha)+\ln(\beta)}\right) \right)}{\ln \alpha}
\end{align*}
and it can be shown that the maximum among the natural numbers $n \in \mathbb{N}$ will occur at the ceiling or floor of $n^{*}$. \qed

\section{Generalizations and Discussion}\label{generalizationN}
In this section, we present a number of generalizations and discussions which will be summarized briefly with the purpose of making the paper more concise and accessible.

\subsection{$N$-step Observability and Its Limitations Due to Policy Dependence} \label{generalizedObsContN}
Given the definition of one step observability in Definition \ref{one_step_observability}, one could conceive of a definition of multiple step observability: for a policy $\gamma$, the POMDP could be called $N$ step observable if for every $f \in C_{b}(\mathcal{X})$ and every $\epsilon>0$ there exists a measurable and bounded function $g$ of $N$ measurements such that
\begin{align*}
\|f(\cdot)-\int_{\mathcal{Y}^N }g(y_{[1,N]})P^{\mu,\gamma}(dy_{[1,N]}|X_{1}=\cdot)\|_{\infty}<\epsilon,
\end{align*}
adopting (\ref{NstepObs}) (note that $N$ can also depend on $f, \epsilon$ if we replace $N$-step observability with observability), which was shown to be an appropriate and consequential observability definition for control-free systems in \cite{mcdonald2018stability}.

However, this approach cannot be used to prove stability due to the past dependency of the control policy. In the following, we explain why the proof method we present is not applicable for such systems unless one restricts the control policies considered. To apply the general process Lemma \ref{general_process_lemma} we need: 
\begin{enumerate}
\item $Y_{[n,n+N-1]}$ can only depend on $Y_{[0,n-1]}$ through $X_{n}$:
\begin{align}
P^{\mu,\gamma}(Y_{[n,n+N-1]}\in \cdot|X_{n},Y_{[0,n-1]})=P^{\mu,\gamma}(Y_{[n,n+N-1]}\in \cdot|X_{n})\label{cont_cond_indp}
\end{align}
\item The channel $Y_{[n,n+N-1]}|X_{n}$ is measure equivalent for $P^{\mu,\gamma}$ and $P^{\nu,\gamma}$:
\begin{align}
P^{\mu,\gamma}(Y_{[n,n+N-1]}\in \cdot|X_{n})&=P^{\nu,\gamma}(Y_{[n,n+N-1]}\in \cdot|X_{n})\label{cont_measure_equiv}
\end{align}
\item The channel $Y_{[n,n+N-1]}|X_{n}$ is time homogeneous:
\begin{align}
P^{\mu,\gamma}(Y_{[n,n+N-1]}\in \cdot|X_{n})&=P^{\mu,\gamma}(Y_{[0,N-1]}\in \cdot|X_{0})~~\forall~n \in \mathbb{N}\label{control_time_homo|}
\end{align}
\end{enumerate}
If we take the LHS of (\ref{cont_cond_indp}) we have
\begin{align*}
&P^{\mu,\gamma}(Y_{[n,n+N-1]}\in \cdot|X_{n},Y_{[0,n-1]}) \\
& =\int_{\mathcal{U}^{N-1}}P^{\mu,\gamma}(Y_{[n,n+N-1]} \in \cdot|X_{n},Y_{[0,n-1]},U_{[n,n+N-2]})P^{\mu,\gamma}(du_{[n,n+N-2]}|X_{n},Y_{[0,n-1]})
\end{align*}
by chain rule of conditional probability. Now it is true $Y_{[n,n+N-1]}|X_{n},U_{[n,n+N-2]}$ is independent of $Y_{[0,n-1]}$ so we can stop conditioning on the past measurements in the inner argument. However, in the outer conditional measure the $U_{[n,n+N-2]}$ may still depend on the past and we have
\begin{align}
\int_{\mathcal{U}^{N-1}}P^{\mu,\gamma}(Y_{[n,n+N-1]} \in \cdot|X_{n},U_{[n,n+N-2]})P^{\mu,\gamma}(du_{[n,n+N-2]}|X_{n},Y_{[0,n-1]})\label{control_past_depand}
\end{align}
If we take the RHS of (\ref{cont_cond_indp}) we have
\begin{align}
P^{\mu,\gamma}(Y_{[n,n+N-1]}\in \cdot|X_{n}) =\int_{\mathcal{U}^{N-1}}P^{\mu,\gamma}(Y_{[n,n+N-1]} \in \cdot|X_{n},U_{[n,n+N-2]})P^{\mu,\gamma}(du_{[n,n+N-2]}|X_{n})\nonumber 
\end{align}
these two equations are not equal for a general control policy, therefore the process fails Definition \ref{general_conditional_indep}, $Y_{[n,n+N-1]}$ does not depend on $Y_{[0,n-1]}$ through $X_{n}$.
As a result, a definition of $N>1$ step observability is incompatible with the proof technique outlined in Lemma \ref{general_process_lemma} and we cannot utilize $N>1$ step observability in a controlled environment to prove stability unless control policies are restricted \sy{e.g. to open-loop control policies where past dependence is avoided, such as sampled control policies where the control is open-loop in between sampling periods, effectively making the measurement a multi-step one. This discussion reveals a further, but not surprising \cite{bartse74}, layer of complexity for the theory of non-linear controlled stochastic systems}.

%

\subsection{Localized Definition of Observability}
As in the control-free case considered in  \cite{mcdonald2018stability}, our observability definition can be generalized further to make its applicability for non-compact domains more general. \sy{While the definition of observability that we introduced is valid for both compact and non-compact state spaces, it is in general difficult to satisfy the definition in a non-compact state space; the concern is that when we have a non-compact space and consider functions $f$ which are not bounded Lipschitz, it may in general be difficult to find a bounded function $g$ which will approximate $f$ through (\ref{defnObsOne}) sufficiently well over the entire state space. }\cm{See \cite[Section 3]{mcdonald2018stability} for some examples of observable and locally observable channels.}


\begin{definition}\label{local_predictableRelaxed}
A POMDP is called {\it locally predictable} (universal in control policies) if there exists a sequence of $\mathcal{F}^{\mathcal{Y}}_{[0,n-1]}$ measurable mappings $a_{n}:\mathcal{Y}^{n}\to \mathcal{X}$ such that, for any $\gamma \in \Gamma$, the family of measures
\begin{align*}
\tilde{\pi}^{\nu,\gamma}_{n-}(\cdot) := \pi_{n-}^{\nu,\gamma}(\cdot+a_{n})
\end{align*}
for every $\mu \ll \nu$, is a uniformly tight family of measures.
\end{definition}

\begin{definition}\label{local_observabilityRelaxed}

A POMDP is called locally observable (universal in control policies) if for every continuous and bounded function $f$, every compact set $K$, every sequence of numbers $a_{n}$, and every $\epsilon>0$, there exists a sequence of uniformly bounded measurable functions $g_{n}$ such that
\sy{
\begin{align*}
&\sup_{x \in K+a_{n}} \left|f(x)-\int_{\mathcal{Y}}g_{n}(y)Q(dy|x) \right|\leq \epsilon, \quad \forall~n \in \mathbb{N} \\
&\sup_{x \not\in K+a_{n}} \left|\int_{\mathcal{Y}}g_{n}(y)Q(dy|x) \right|\leq 2\|f\|_{\infty}, \quad 
\forall~ n \in \mathbb{N}
\end{align*}
}

\end{definition}

Then, we can state the following generalization of Theorem \ref{weak_merging_pred} (whose proof follows similarly).
\begin{theorem}\label{local_thmRelaxed}
Assume $\mu \ll \nu$ and that the POMDP is locally predictable and locally observable (universal in control policies). Then the predictor is universally stable weakly a.s. .
\end{theorem}


\subsection{Robustness Under Weak Merging of Priors}
Our robustness result is also applicable under weak convergence of priors, provided that the channel is continuous under total variation, and the filter is stable subject to a modest rate of convergence condition. These would build on \cite[Theorem 3.3]{kara2019robustness} and the analysis presented in Section \ref{robustMainSec}.

\section{Conclusion}\label{sec:conc}
Filter stability has been studied extensively in control-free contexts.  In this paper, we studied the filter stability problem, developed new methods and results for controlled stochastic dynamical systems, and studied the implications of filter stability on robustness of optimal solutions for partially observed stochastic control problems. 
\section{Appendix}\label{Appendix_Observability_Predictor}
\subsection{Proof of Theorem \ref{avg_cost_robust}}
Consider the robustness difference
\begin{align}
J_{\infty}(\mu,\gamma^{\nu})-J_{\infty}^{*}(\mu)
&=J_{\infty}^{*}(\nu)-J_{\infty}^{*}(\mu)+J_{\infty}(\mu,\gamma^{\nu})-J_{\infty}^{*}(\nu)\nonumber\\
&\leq |J_{\infty}^{*}(\nu)-J_{\infty}^{*}(\mu)|+J_{\infty}(\mu,\gamma^{\nu})-J_{\infty}(\nu,\gamma^{\nu})\label{continuity_policy_match}
\end{align}
\cm{Consider now the difference $|J_{\infty}^{*}(\mu)-J_{\infty}^{*}(\nu)|$. If $J_{\infty}^{*}(\mu)<J_{\infty}^{*}(\nu)$ then  the absolute value of their difference is the larger value $J_{\infty}^{*}(\nu)$ minus the smaller value $J_{\infty}^{*}(\mu)$. Then since $J_{\infty}^{*}(\nu)=J_{\infty}(\nu,\gamma^{\nu})$, if we replace $\gamma^{\nu}$ with $\gamma^{\mu}$ we have $J_{\infty}(\nu,\gamma^{\mu})>J_{\infty}^{*}(\nu)$ and the difference is even greater. Therefore when $J^{*}_{\infty}(\mu)<J^{*}_{\infty}(\nu)$ we have,
\begin{align*}
|J_{\infty}^{*}(\mu)-J_{\infty}^{*}(\nu)|&= J^{*}_{\infty}(\nu)-J^{*}_{\infty}(\mu)\\
&=J_{\infty}(\nu, \gamma^{\nu})-J_{\infty}(\mu,  \gamma^{\mu})\\
&\leq J_{\infty}(\nu,\gamma^{\mu})-J_{\infty}(\mu,\gamma^{\mu})
\end{align*}}
By a symmetric argument for the other direction in the inequality, we arrive at
\begin{align*}
|J_{\infty}^{*}(\mu)-J_{\infty}^{*}(\nu)| \leq \max\left(J_{\infty}(\mu,\gamma^{\nu})-J_{\infty}(\nu,\gamma^{\nu}),J_{\infty}(\nu,\gamma^{\mu})-J_{\infty}(\mu,\gamma^{\mu}) \right)
\end{align*}
and we can ultimately determine the robustness bound by studying the expected cost operator $J_{\infty}(\cdot,\gamma^{\nu})$ under different priors but the same control policy. By same control policy, we mean $\gamma^{\nu}$ is the optimal control policy designed with respect to the prior $\nu$. This means the DM sees observations $y_{[0,n]}$, computes the filter believing the prior is $\nu$ (the mapping from measurements to control actions is the same for each POMDP with different priors). 

Note that for two non-negative bounded sequences $0<a_{n}<m<\infty$ and $0<b_{n}<k<\infty$ we have that the difference of their limsup's is less than the limsup of the difference $\limsup_{n \to \infty} a_{n}-\limsup_{n \to \infty}b_{n}\leq \limsup_{n \to \infty}(a_{n}-b_{n})$, 
therefore since $c$ is non-negative and bounded we have
\begin{align*}
&J_{\infty}(\mu,\gamma^{\nu})-J_{\infty}(\nu,\gamma^{\nu})\\
=&\limsup_{T \to \infty}\frac{1}{T}\left(\sum_{i=0}^{T-1}E^{\mu,\gamma^{\nu}}[c(X_{i},U_{i})] \right) -\limsup_{T \to \infty}\frac{1}{T}\left(\sum_{i=0}^{T-1}E^{\nu,\gamma^{\nu}}[c(X_{i},U_{i})] \right)\\
\leq&\limsup_{T \to \infty}\frac{1}{T}\sum_{i=0}^{T-1}\left(E^{\mu,\gamma^{\nu}}[c(X_{i},U_{i})]-E^{\nu,\gamma^{\nu}}[c(X_{i},U_{i})] \right) \\
&\leq \limsup_{T \to \infty}\frac{\|c\|_{\infty}}{T}\left(\sum_{i=0}^{T-1}\|P^{\mu,\gamma^{\nu}}((X_{i},U_{i}) \in \cdot)  -P^{\nu,\gamma^{\nu}}((X_{i},U_{i})\in \cdot)\|_{TV} \right)
\end{align*}
then we see that
\begin{small}
\begin{align}
&\|P^{\mu,\gamma^{\nu}}((X_{i},U_{i}) \in \cdot)-P^{\nu,\gamma^{\nu}}((X_{i},U_{i})\in \cdot)\|_{TV} \nonumber\\
=&
\sup_{\|f\|_{\infty}\leq 1} \left|\int_{\mathcal{X} \times \mathcal{U}}f(x,u)P^{\mu,\gamma^{\nu}}(dx_{i},du_{i})-\int_{\mathcal{X} \times \mathcal{U}}f(x,u)P^{\nu,\gamma^{\nu}}(dx_{i},du_{i})\right| \nonumber\\
=&\sup_{\|f\|_{\infty}\leq 1} \left|\int_{\mathcal{X}}\int_{\mathcal{X} \times \mathcal{U}}f(x,u)P^{\mu,\gamma^{\nu}}(dx_{i},du_{i}|X_{0})\mu(dx_{0})  -\int_{\mathcal{X}}\int_{\mathcal{X} \times \mathcal{U}}f(x,u)P^{\nu,\gamma^{\nu}}(dx_{i},du_{i}|X_{0})\nu(dx_{0})\right| 
\end{align}
As was discussed, both POMDP use the same control policy, which maps measurements to control actions in the same way. Once we fix the realization of $X_{0}=x$ in either case, the distribution on $Y_{0}$ is the same, hence the distribution on $U_{0}$ is the same, and hence $X_{1},Y_{1},U_{1}$ and so on. Therefore in the above the two inner integrals are the same function of $x$ and we can upper bound by $\|\mu-\nu\|_{TV}$.
\end{small}
Thus, $J_{\infty}(\mu,\gamma^{\nu})-J_{\infty}(\nu,\gamma^{\nu}) \leq \|c\|_{\infty}\|\mu-\nu\|_{TV}$.
Both terms in (\ref{continuity_policy_match}) have this bound and the overall bound is multiplied by a factor of 2.
\qed

\subsection{Proof of Lemma \ref{tv_structure}}~
\begin{lemma}\label{restrict_structure}\cite[Lemma 4.6]{mcdonald2018stability}
Assume $\mu \ll \nu$ and fix any control policy $\gamma$. For any sigma field $\mathcal{G} \subseteq \mathcal{F}^{\mathcal{X}}_{0,\infty} \vee \mathcal{F}^{\mathcal{Y}}_{0,\infty}$ we have:
\begin{align*}
\frac{dP^{\mu,\gamma}|_{\mathcal{G}}}{dP^{\nu,\gamma}|_{\mathcal{G}}}&=E^{\nu,\gamma}\left[\left.\frac{d\mu}{d\nu}(X_{0})\right|\mathcal{G}\right]~~~P^{\mu,\gamma}~a.s.
\end{align*}
\end{lemma}
\begin{proof}
\cm{Note that conditioned on the value of $X_{0}$, knowledge of the prior is irrelevant so for any set $A \in \mathcal{G}$ we have
\begin{align*}
E^{\mu}[1_{A}] = E^{\nu}[1_{A}\frac{d\mu}{d\nu}(X_{0})]
\end{align*}
then we can apply law iterated expectations, condition internally on $\mathcal{G}$ and move the indicator out in the inner expectation since it is $\mathcal{G}$ measurable. This shows
\begin{align*}
E^{\mu}[1_{A}]&=E^{\nu}[1_{A}E^{\nu}[\frac{d\mu}{d\nu}(X_{0})|\mathcal{G}]]
\end{align*}}
\end{proof}

\begin{lemma}\label{condition_structure}\cite[Lemma 4.7]{mcdonald2018stability}
Assume $\mu \ll \nu$ and fix any control policy $\gamma$. For any two sigma fields $\mathcal{G}_{1},\mathcal{G}_{2} \subset \mathcal{F}^{\mathcal{X}}_{0,\infty} \vee \mathcal{F}^{\mathcal{Y}}_{0,\infty}$, let $P^{\mu,\gamma}|_{\mathcal{G}_{1}}|\mathcal{G}_{2}$ represent the probability measure $P^{\mu,\gamma}$ restricted to $\mathcal{G}_{1}$, conditioned on field $\mathcal{G}_{2}$. We then have
\begin{align*}
\frac{dP^{\mu,\gamma}|_{\mathcal{G}_{1}}|\mathcal{G}_{2}}{dP^{\nu,\gamma}|_{\mathcal{G}_{1}}|\mathcal{G}_{2}}&=\frac{E^{\nu,\gamma}[\frac{d\mu}{d\nu}(X_{0})|\mathcal{G}_{1} \vee \mathcal{G}_{2}]}{E^{\nu,\gamma}[\frac{d\mu}{d\nu}(X_{0})|\mathcal{G}_{2}]}~~~P^{\mu,\gamma}~a.s.
\end{align*}
\end{lemma}

\begin{proof}
\cm{ 
Note for any set $A \in \mathcal{G}_{1}$ we can write $E^{\mu}[1_{A}]$ in two ways
\begin{align*}
E^{\mu}[1_{A}]&=E^{\nu}[E^{\mu}[1_{A}\frac{dP^{\mu}|_{\mathcal{G}_{2}}}{dP^{\nu}|_{\mathcal{G}_{2}}}|\mathcal{G}_{2}]]\\
E^{\mu}[1_{A}]&=E^{\nu}[1_{A}E^{\nu}[\frac{dP^{\mu}}{dP^{\nu}}|\mathcal{G}_{1}\vee \mathcal{G}_{2}]]
\end{align*}
applying Lemma 3.9
\begin{align*}
E^{\mu}[1_{A}]&=E^{\nu}[E^{\mu}[1_{A}E^{\nu}[\frac{d\mu}{d\nu}(X_{0})|\mathcal{G}_{2}]|\mathcal{G}_{2}]]\\
E^{\mu}[1_{A}]&=E^{\nu}[E^{\nu}[1_{A}\frac{d\mu}{d\nu}(X_{0})|\mathcal{G}_{1}\vee \mathcal{G}_{2}]]
\end{align*}
by definition the Radon Nikodym derivative $\frac{dP^{\mu}|_{\mathcal{G}_{1}}|\mathcal{G}_{2}}{dP^{\nu}|_{\mathcal{G}_{1}}|\mathcal{G}_{2}} = f$ is the function that satisfies

\begin{align*}
E^{\nu}[E^{\nu}[1_{A}fE^{\nu}[\frac{d\mu}{d\nu}(X_{0})|\mathcal{G}_{2}]|\mathcal{G}_{2}]]&=E^{\mu}[1_{A}]=E^{\nu}[E^{\nu}[1_{A}\frac{d\mu}{d\nu}(X_{0})|\mathcal{G}_{1}\vee \mathcal{G}_{2}]]
\end{align*}
one can show substituting $f$ for the function in the theorem has this property.
}
\end{proof}

\begin{proof}[Proof of Lemma \ref{tv_structure}]
\cm{ 
For to measures $\mu_{1}$ and $\mu_{2}$ note total variation can be written as
\begin{align*}
\|\mu_{1}-\mu_{2}\|_{TV} = E^{\mu_{2}}\left[\left|\frac{d\mu_{1}}{d\mu_{2}}-1 \right|\right]
\end{align*}
Applying Lemma 3.10 and cross multiplying to get a single fraction in the absolute value leads to the result.
}
\end{proof}



\end{document}